%% file: 0-main.tex
\documentclass[acmtog,screen,nonacm]{acmart}

\usepackage{booktabs} % For formal tables
\setcitestyle{square}
\usepackage{ifthen}
\usepackage{enumitem}
\usepackage{algorithm}
\usepackage[noend]{algpseudocode}
\usepackage{syntax}
\usepackage{amsfonts}
\usepackage{listings}
\usepackage{fancyvrb}
\usepackage{wrapfig}
\usepackage{graphicx}
\usepackage{xspace}
\usepackage{colortbl}
\usepackage{gensymb}
\usepackage{verbatim}
\usepackage{colortbl}
\usepackage{booktabs}
\usepackage{multirow}
\usepackage{cleveref}
\usepackage{pifont}
\usepackage{circledsteps}
\usepackage{siunitx}

\usepackage{makecell}  % 引入makecell包

\usepackage{nomencl}
\makenomenclature

\algdef{SE}[DOWHILE]{Do}{doWhile}{\algorithmicdo}[1]{\algorithmicwhile\ #1}%

\definecolor{mygreen}{rgb}{0,0.6,0}                         
\definecolor{mygray}{rgb}{0.95,0.95,0.95}
\definecolor{codebg}{rgb}{0.95, 0.95, 0.95}  
\newcommand{\wz}[1]{{\color{black} #1}}

\usepackage{amsthm,amsmath,amssymb}
\usepackage{mathrsfs}
\copyrightyear{2025}
\acmYear{2025}
\setcopyright{acmlicensed}
\acmConference[SIGGRAPH Conference Papers '25]{Special Interest Group on Computer Graphics and Interactive Techniques Conference Conference Papers }{August 10--14, 2025}{Vancouver, BC, Canada}
\acmBooktitle{Special Interest Group on Computer Graphics and Interactive Techniques Conference Conference Papers (SIGGRAPH Conference Papers '25), August 10--14, 2025, Vancouver, BC, Canada}
\acmDOI{10.1145/3721238.3730700}
\acmISBN{979-8-4007-1540-2/2025/08}

\begin{document}

\title{DualMS: Implicit Dual-Channel Minimal Surface Optimization for Heat Exchanger Design}

\author{Weizheng Zhang}
\email{wonderz.top@gmail.com}
\affiliation{
  \institution{Shandong University}
  \country{China}
}
\author{Hao Pan}
\email{haopan@tsinghua.edu.cn}
\affiliation{
  \institution{Tsinghua University}
  \country{China}
}
\author{Lin Lu}
\authornote{corresponding author}
\email{llu@sdu.edu.cn}
\affiliation{
  \institution{Shandong University}
  \country{China}
}
\author{Xiaowei Duan}
\email{duanxwss@gmail.com}
\affiliation{
  \institution{Shandong University}
  \country{China}
}
\author{Xin Yan}
\email{yanxin@sdu.edu.cn}
\affiliation{
  \institution{Shandong University}
  \country{China}
}
\author{Ruonan Wang}
\email{wangruonan@iet.cn}
\affiliation{
  \institution{Institute of Engineering Thermophysics, Chinese Academy of Sciences}
  \country{China}
}
\author{Qiang Du}
\email{duqiang@iet.cn}
\affiliation{
  \institution{Institute of Engineering Thermophysics, Chinese Academy of Sciences}
  \country{China}
}

\begin{abstract}

Heat exchangers are critical components in a wide range of engineering applications, from energy systems to chemical processing, where efficient thermal management is essential. The design objectives for heat exchangers include maximizing the heat exchange rate while minimizing the pressure drop, requiring both a large interface area and a smooth internal structure. State-of-the-art designs, such as triply periodic minimal surfaces (TPMS), have proven effective in optimizing heat exchange efficiency. However, TPMS designs are constrained by predefined mathematical equations, limiting their adaptability to freeform boundary shapes. Additionally, TPMS structures do not inherently control flow directions, which can lead to flow stagnation and undesirable pressure drops.

This paper presents \emph{DualMS}, a novel computational framework for optimizing dual-channel minimal surfaces specifically for heat exchanger designs in freeform shapes. To the best of our knowledge, this is the first attempt to directly optimize minimal surfaces for two-fluid heat exchangers, rather than relying on TPMS. Our approach formulates the heat exchange maximization problem as a constrained connected maximum cut problem on a graph, with flow constraints guiding the optimization process. To address undesirable pressure drops, we model the minimal surface as a classification boundary separating the two fluids, incorporating an additional regularization term for area minimization. We employ a neural network that maps spatial points to binary flow types, enabling it to classify flow skeletons and automatically determine the surface boundary. 
DualMS demonstrates greater flexibility in surface topology compared to TPMS and achieves superior thermal performance, with lower pressure drops while maintaining a similar heat exchange rate under the same material cost.
The project is open-sourced at \href{https://github.com/weizheng-zhang/DualMS}{https://github.com/weizheng-zhang/DualMS}.

\end{abstract}



\begin{teaserfigure}
\includegraphics[width=1.0\linewidth]{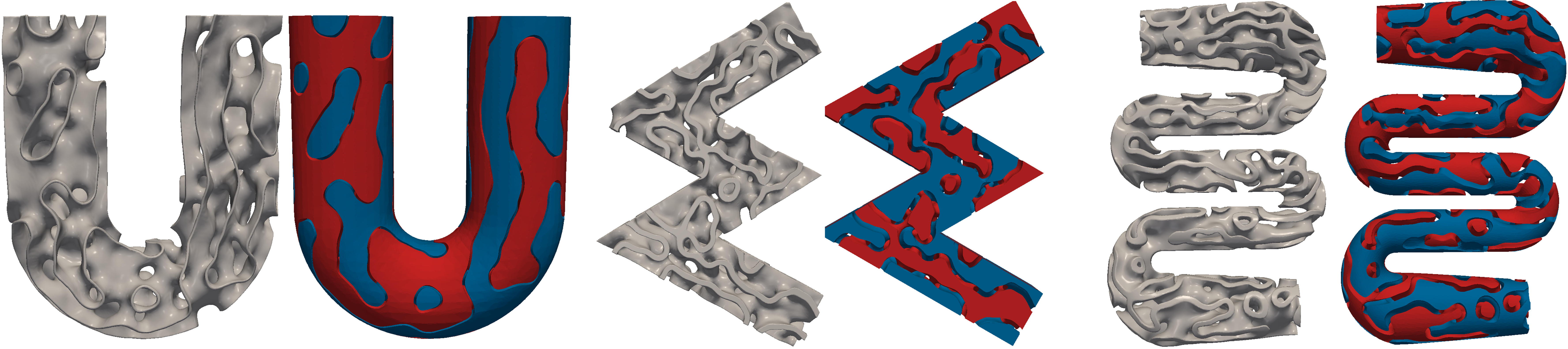}
\caption{DualMS generates an optimized separation surface (left) to enhance heat exchange between two fluids (right) within freeform domains.}
\label{fig:teaser}    
\end{teaserfigure}

\maketitle

\input{1-intro}
\input{2-related}
\input{3-overview}
\input{4-method}
\input{5-results}
\input{6-conclusion}

\begin{acks}
We thank all reviewers for their valuable comments and constructive suggestions. The authors acknowledge the support of the National Natural Science Foundation of China (NSFC) through the Excellence Research Group Program (Grant No.~52488101). This work was also supported in part by the National Natural Science Foundation of China (Grant Nos.~62472258 and 62302275).
\end{acks}

\bibliographystyle{formats/ACM-Reference-Format}
\bibliography{DualMS} 

\newpage

% \maketitle

In the supplementary materials, we provide detailed information on computational performance, the convergence of the maximum cut algorithm, numerical simulation details, an ablation study on the flow field, and fabrication results.

\input{appendix}

\nomenclature{\(T\)}{Temperature}
\nomenclature{\(P\)}{Pressure}
\nomenclature{\(\dot{m}\)}{Mass Flow Rate}
\nomenclature{\(\rho\)}{Density}
\nomenclature{\(Re\)}{Reynolds Number, \(Re=\frac{\rho_0 u L}{\mu_0}\)}
\nomenclature{\(\kappa\)}{Thermal Conductivity}
\nomenclature{\(\mu\)}{Dynamic Viscosity}
\nomenclature{\(C_p\)}{Specific Heat Capacity}
\nomenclature{\(L\)}{Length}
\nomenclature{\(_{in}\)}{Value at the Inlet}
\nomenclature{\(_{out}\)}{Value at the Outlet}
\nomenclature{\(_0\)}{Reference Value}
\nomenclature{\(t\)}{Time}
\nomenclature{\(\tau\)}{Shear Stress}
\nomenclature{\(e\)}{Internal Energy}

\end{document}

%% file: 1-intro.tex
\section{Introduction}
\label{sec:intro}

Heat exchangers are essential devices for transferring heat between two or more fluids, typically separated by solid walls to prevent mixing. They are widely used in aerospace, automotive, and industrial applications, where efficient thermal management is crucial for optimizing performance, ensuring reliability, and improving energy efficiency~\cite{Zohuri2017book}. The primary design objective of heat exchangers is to maximize heat transfer while minimizing flow resistance, achieving an optimal balance between thermal efficiency and pressure drop.

The advent of additive manufacturing (AM) has revolutionized heat exchanger design by enabling the fabrication of complex structures that were previously impossible with traditional methods. This has also led to the adoption of freeform shapes in heat exchangers, enabling higher integration in spaces with complex, non-rectilinear geometries, such as those found in aerospace or automotive systems, where conventional designs cannot efficiently utilize available volume or adapt to such spaces.
 
Among the advancements in heat exchanger design, triply periodic minimal surfaces (TPMS) have gained significant attention due to their high surface-area-to-volume ratios and well-distributed internal channels, which enhance heat transfer and optimize fluid flow dynamics \cite{Oh2023, Gado2024}.
These advantages make TPMS designs superior to conventional types, such as plate, plate-fin, and printed circuit designs \cite{Ayub2003, Li2011, Ngo2006}. 
However, TPMS structures are typically constrained by their regular spatial tiling patterns and limited adaptability to flow directions, restricting their applicability in freeform domains. 

To address these limitations, this paper introduces \emph{DualMS}, an implicit dual-channel minimal surface optimization framework tailored for heat exchanger design in freeform shapes. The key idea behind DualMS is to directly optimize the minimal surfaces that separate the two fluids, bypassing the reliance on TPMSs and thereby achieving a greater degree of freedom in the channel geometries.

The objective of \emph{DualMS} is twofold: achieving high heat exchange efficiency and minimizing flow resistance. High heat exchange efficiency requires a large separation surface area between the two flow channels, while low flow resistance demands a smooth and continuous minimal surface. 
To balance these inherently conflicting objectives, we propose a two-stage optimization process.
In the first stage, we globally maximize the surface area by performing a rough spatial partition to enhance heat transfer efficiency. 
In the second stage, we minimize the surface area locally to reduce pressure drop, subject to skeleton constraints obtained by the global area maximization, to reduce flow resistance.

To globally maximize the surface area separating two flows, we first generate a dual-flow skeleton designating the overall interaction  between the two flow channels. 
Specifically, we construct a graph within the design space and formulate it as a connected maximum cut problem with flow constraints. This approach partitions the graph by cutting edges, optimizing the interconnections between the two resulting subgraphs, where each subgraph serves as the flow skeleton for one of the fluids.

For the local minimization of pressure drop, we aim to use minimal surfaces as the structural geometry while treating the dual-flow skeleton as hard constraints. These surfaces are designed to be smooth and simple, reducing turbulence and vortices in the fluid flow, thereby improving the flow dynamics.
Unlike the classical Plateau's problem with fixed boundaries, our task requires solving for a minimal surface that separates two flow channels with an unknown boundary. To address this, we model the minimal surface as the classification boundary between the flow channels, incorporating a regularization term to minimize surface area. A neural network is trained to classify spatial points based on the flow skeletons while minimizing the decision boundary area, enabling it to adapt to complex geometries and automatically determine the surface boundary.

Our main contributions are as follows:  
\begin{itemize}[nosep]
    \item We approach the heat exchanger design problem by directly optimizing the separation surface of two fluids, focusing on achieving the optimal balance between surface area maximization and pressure drop minimization.
    
    \item We frame the surface optimization as a two-stage strategy: first, globally determining the flow topology of the two fluids, and then locally  optimizing the interface minimal surface as the decision boundary of a binary classifier model.
    
    \item We propose \emph{DualMS}, a novel computational framework for dual-channel heat exchanger design with freeform boundary shapes. This framework offers greater flexibility in surface topology compared to TPMS (Figure~\ref{fig:teaser}) and demonstrates superior thermal performance by achieving a lower pressure drop while maintaining a comparable heat exchange rate under the same material cost.

\end{itemize}  

%% file: 2-related.tex
\section{Related Work}
\label{sec:related}

\begin{figure*}[t]
    \centering
    \includegraphics[width=1\linewidth]{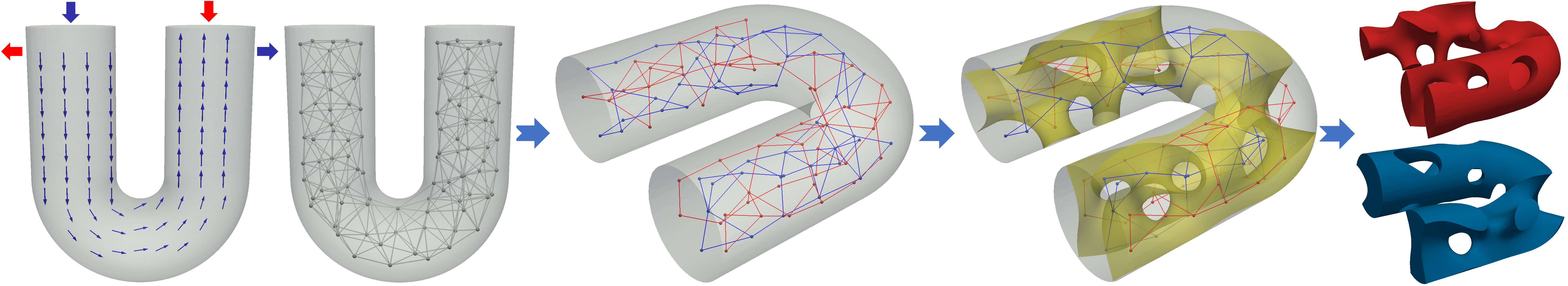}
        \leftline{ \footnotesize  \hspace{0.07\linewidth}
            (a)  \hspace{0.15\linewidth} 
            (b) \hspace{0.19\linewidth}
            (c)  \hspace{0.25\linewidth}
            (d) \hspace{0.17\linewidth}
            (e)}
        \caption{Given the design domain and boundary conditions and  of the heat exchanger (a), DualMS initializes an undirected graph (b) derived from centroidal Voronoi tessellations and the given flow field. The method then optimizes the dual flow skeletons by solving a constrained connected maximum cut problem on the graph (c), and then  optimizes a minimal surface (d) to effectively partition the domain into two distinct flow channels (e).}
    \label{fig:overview}
\end{figure*}

\paragraph{Minimal Surface Modeling}

Minimal surfaces are variational surfaces with locally minimal surface areas. As a long-studied subject, they are both theoretically intriguing and practically important.
Initially formulated by Joseph Plateau~\cite{plateau1873statique}, minimal surfaces are assumed naturally by soap films spanning curved wire boundaries due to surface tension.

To compute minimal surfaces, one can use known constructions of fundamental minimal patches combined through inherent symmetries or apply numerical schemes to solve discrete approximations. 
The former leads to the construction of TPMS \cite{1970Gesammelte,Schoen1970,FischerKoch1987,karcher1996construction}, while the latter typically uses a mesh-based surface representation driven by mean curvature flow to produce minimal surface meshes \cite{Pinkall1993,Brakke1992}.

TPMS are typically constrained by their regular spatial tiling patterns, which limits their adaptability to freeform domains. Nevertheless, they are widely used in applications such as porous materials and thermal management systems due to their local smoothness, high surface-area-to-volume ratio, and superior overall mechanical properties~\cite{Liu2021}. Researchers have also attempted to optimize the periodicity of implicit functions to better conform to freeform shapes, enabling enhanced structural performance or energy absorption~\cite{Yan2019,Tian2024}. Recent works have further expanded the design space of TPMS by exploring the boundary skeletons of surface patches and applying corresponding symmetry rules~\cite{Makatura2023,Xu2023}.
In this paper, we focus on designing more flexible minimal surfaces without strict regularity constraint, to adapt to different domain shapes.

The numerical solving of minimal surfaces represented by discrete meshes can be achieved via finite difference methods \cite{Hinata1974} or the mean curvature flow \cite{Brakke1992, Desbrun1999, Xu2006}.
However, the mesh quality during the surface evolution can deteriorate significantly, requiring frequent remeshing \cite{Brakke1992} or modified energies that are mesh quality aware \cite{Pan2012}.
Moreover, challenges like free boundaries, where the boundary curve of an evolved surface can change along a given domain, as well as changing surface topologies, make the mesh-based evolution even more challenging \cite{Brakke1992}.

\cite{Wang2021} proposes to solve minimal surfaces via the framework of geometric measure theory (GMT), where curves and surfaces are represented implicitly as currents that generalize to be volumetric distributions. 
In this setting, minimal surfaces become the solutions of convex minimal mass norm problems, and can be solved over the volumetric domain, bypassing mesh qualities or changing topologies.
DeepCurrents \cite{palmer2022deepcurrents} further takes the idea into a neural network based current representation, enabling high-quality and more robust computation.
In this paper, we also use neural network based implicit fields to model minimal surfaces, to achieve high surface quality and robustness to changing boundaries and topology.
However, we differ from these previous works by solving a different kind of input than given boundary curves as specified by the classical Plateau's problem.
Instead, given the heat exchanger design problem, we seek minimal surfaces separating two different types of flow channels, without knowing the boundary curves.
We solve this problem by formulating the minimal surfaces as the decision boundary of a binary classifier model, and apply a total variation loss to ensure the decision boundary is locally area minimizing, which is theoretically grounded in GMT.

Related to minimal surface modeling, \cite{Palmer2024} solves the problem of direction field synthesis by formulating it as a minimal surface problem in the lifted space of circle bundles over a given domain.
\cite{Perez2017} solves a variant of the Plateau's problem, where the minimal surface tensile force is balanced by a given elastic rod.

\paragraph{Heat Exchanger Design}

Early efforts in heat exchanger design using topology optimization (TO) include the works of \cite{Hoeghoej2020,Kobayashi2020}, which employed a single design variable field to optimize two-fluid heat exchangers. While these studies showed promising results, they highlighted a key limitation of the TO framework: the difficulty in enforcing manufacturing constraints, such as minimal feature size, often resulting in impractical designs for physical fabrication \cite{Fawaz2022}. 
Additionally, many flow optimization studies have explored various types of flows, such as Stokes flow~\cite{Borrvall2003,Li2022} and viscous flow~\cite{Kontoleontos2013}, highlighting the complexity of integrating flow dynamics into heat exchanger design.

Further advancements were made by \cite{Feppon2021}, who employed the level set method to optimize 2D and 3D heat exchangers involving two fluids. This approach allowed for more flexibility in design but still faced limitations related to managing complex geometries and flow distributions.

Recent studies have increasingly explored the use of TPMS for heat exchangers, taking advantage of their ability to partition 3D space into smooth, continuous, non-intersecting domains. Most research focuses on analyzing the heat transfer and pressure drop characteristics of TPMS structures, considering factors such as TPMS type and wall thickness \cite{Kim2020,Kaur2021, Attarzadeh2021, Iyer2022, Reynolds2023,Wang2024}. 
Some efforts have extended this by optimizing TPMS parameters, including periodicity and wall thickness \cite{Jiang2023}, though these approaches remain confined to the degrees of freedom allowed by TPMS. Despite their advantages, TPMS-based heat exchangers are generally limited to unidirectional fluid flow configurations. To improve this, \cite{Oh2023} introduced three filters to enhance TPMS adaptability: a selection filter for inlets and outlets, a barrier filter to control flow direction, and a boundary filter to reduce flow resistance. However, the local frictional losses around the barriers can increase flow energy dissipation.
In contrast, our approach optimizes surfaces directly based on flow field guidance, offering greater flexibility in surface topology compared to TPMS structures, and achieving lower pressure drop while maintaining a similar heat exchange rate.

%% file: 3-overview.tex
\section{Overview}
\label{sec:Overview}

Given a freeform heat exchanger design domain with specified boundary conditions (inlet and outlet positions) and a predefined flow field providing rough flow guidance (Figure~\ref{fig:overview}a), our objective is to design internal surfaces that balance high heat exchange efficiency and low flow resistance. High heat exchange efficiency corresponds to a large separation surface area between the two flow channels, while low flow resistance requires a smooth, continuous minimal surface.

To address this dual-objective max-min optimization problem, we propose a two-stage process: \emph{dual flow skeleton optimization} and \emph{minimal surface optimization}.

Dual flow skeleton optimization seeks to determine the flow topology for the two fluids (Section~\ref{sec:dualSkeleton}). We introduce the concept of a flow skeleton, a graph embedded in the design space that represents fluid flow paths. The heat exchanger’s two fluid channels are modeled by two such flow skeletons, which also define the interface surface separating the two channels. To compute the skeletons, we initialize a dense graph within the design domain (Figure~\ref{fig:overview}b) and solve it as a connected maximum cut problem with flow constraints. This partitions the graph into two subgraphs, each representing one fluid channel. The flow constraints ensure no vertex with degree 1 exists in either subgraph, preventing dead-end paths in the flow channels (Figure~\ref{fig:overview}c).

Minimal surface optimization models the separation surface as a classification boundary with an additional regularization term to minimize surface area (Section~\ref{sec:MSOptimization}). Starting with the initial fluid classification, we train a neural network to classify the design space, aligning with the flow skeletons while minimizing the decision boundary area to form the optimal separation surface (Figure~\ref{fig:overview}d).

Finally, we apply a uniform thickness to achieve a smooth, continuous separating wall, effectively dividing the design space into two distinct flow channels (Figure~\ref{fig:overview}e).

%% file: 4-method.tex
\section{Dual Flow Skeleton Optimization}
\label{sec:dualSkeleton}

In this section, we outline the process of converting the input freeform shape and boundary conditions into an optimized dual-flow skeleton. This involves initializing a graph within the shape domain, assigning edge weights based on the flow field, and solving a constrained connected maximum cut problem to generate the dual-flow skeletons.

\subsection{Graph Initialization}

We initialize an undirected graph \( G(V, E) \) to represent the design space, where \( V \) consists of stochastically sampled points optimized via centroidal Voronoi tessellation (CVT)~\cite{Liu2009}. Edges \( E \) are generated using Delaunay triangulation, with weights assigned to incorporate the flow field by penalizing edges that are more perpendicular to the flow direction.
The number of vertices impacts computational time and the degree of freedom: a higher vertex count provides finer spatial partitioning and greater surface area, enhancing heat transfer but increasing pressure drop and volume fraction, as shown in the results.

For each edge \( (u, v) \in E \), the weight \( w_{uv} \) is computed based on the alignment of the edge direction \( \mathbf{d}_{uv} \) and the average flow vector \( \mathbf{f}_{uv} \). Let \( \theta_{uv} \) be the angle between \( \mathbf{d}_{uv} \) and \( \mathbf{f}_{uv} \). The weight is defined as:
\begin{equation}
    w_{uv} =
    \begin{cases}
        a, & \text{if } \theta_{uv} \in (\frac{\pi}{4}, \frac{3\pi}{4}),\\
        1, & \text{otherwise},
    \end{cases}
\end{equation}
where \( a > 1 \) is a penalty factor, empirically set to 5 in our experiments. This weighting scheme encourages the graph to preserve edges aligned with the flow, guiding fluid channels along the flow direction and separation surfaces perpendicular to it.

\subsection{Dual Skeleton from Connected Maximum Cut}

To maximize the contact surface between the two flow channels, we partition \( G(V, E) \) into two connected subgraphs, \( G_1(V_1, E_1) \) and \( G_2(V_2, E_2) \), representing the two fluid channels of the heat exchanger. The objective is to intertwine the edges of \( G_1 \) and \( G_2 \) as much as possible, i.e., to maximize the sum of the edge weights of the edges removed during partitioning. This is a variant of the \emph{connected maximum cut problem}, as formulated in \cite{haglin1991approximation}.

In the connected maximum cut problem, given an undirected graph \( G = (V, E) \) with edge weights \( w_{uv} \geq 0 \) for \( (u, v) \in E \), the goal is to partition the vertex set \( V \) into two disjoint subsets \( V_1 \) and \( V_2 \) such that the subgraphs \( G_1(V_1, E_1) \) and \( G_2(V_2, E_2) \) are connected, and the total weight of edges crossing the partition is maximized. 

In our formulation, we extend the classical problem by adding a constraint: ensuring that each vertex in both \( G_1 \) and \( G_2 \) has a degree greater than 1, thereby preventing dead-end paths in the flow channels. Let \( x_i \in \{-1, +1\} \) indicate the subset membership of vertex \( i \), where \( x_i = +1 \) if \( i \in V_1 \) and \( x_i = -1 \) if \( i \in V_2 \). The optimization problem is formulated as:
\begin{equation}
\begin{aligned}
\max \quad & \sum_{(u, v) \in E} w_{uv} \cdot \frac{1 - x_u x_v}{2}, \\
\text{subject to} \quad & V_1 \cup V_2 = V, \quad V_1 \cap V_2 = \emptyset, \\
& \text{$G_1(V_1, E_1)$ and $G_2(V_2, E_2)$ are connected,} \\
& \deg_{G_1}(v) > 1 \quad \forall\, v \in V_1,\quad\\
& \deg_{G_2}(v) > 1 \quad \forall\, v \in V_2.
\end{aligned}
\end{equation}

Since the connected maximum cut problem is NP-hard, the constrained version does not have a polynomial-time solution. Therefore, we employ a heuristic algorithm.

The objective function is defined as:
\begin{equation}
        f(G, G_1, G_2) = \sum_{(u, v) \in E} w_{uv} - \left(\sum_{(u, v) \in E_1} w_{uv} + \sum_{(u, v) \in E_2} w_{uv}\right).
\label{eq:maxcut}
\end{equation}
This function represents the sum of the cut edges, which we aim to maximize.

The algorithm proceeds by initially partitioning \( G(V, E) \) into two connected subgraphs, \( G_1(V_1, E_1) \) and \( G_2(V_2, E_2) \). We then iterate through each vertex, considering its move to the opposite subgraph and recording the change in the objective function. The vertex that maximizes this change is swapped, and the process repeats until the objective function converges.

The algorithm can accommodate additional constraints. For instance, we ensure subgraph connectivity by checking whether a vertex swap would disconnect either subgraph. To avoid dead ends, we also check that no vertex in either subgraph has degree 1 after the swap.

\begin{algorithm}[H]
\caption{Constrained Connected Maximum Cut Algorithm}
\label{alg:heuristicMaxCut}
\begin{algorithmic}[1]
\State \textbf{Input:} Graph \( G(V, E) \)
\State \textbf{Output:} Approximate optimal partition \( G_1, G_2 \)

\State Initialize two empty subgraphs \( G_1(V_1, E_1) \) and \( G_2(V_2, E_2) \)
\State Assign vertices in \( G \) to \( G_1 \) and \( G_2 \), ensuring that no connected subgraph contains vertices with degree 1.

\State Compute the initial objective function \( f(G, G_1, G_2) \) with Eq.~\ref{eq:maxcut}

\Repeat
    \For{each vertex \( v \in V \)}
        \State Simulate moving \( v \) to the other subgraph
        \State Check constraints:
        \State \quad Ensure connectivity of \( G_1 \) and \( G_2 \)
        \State \quad Avoid introducing degree-1 vertices
        \State Evaluate the change in \( f(G, G_1, G_2) \) after the swap
    \EndFor
    \State Identify the vertex \( v^* \) whose swap maximizes the increase in \( f(G, G_1, G_2) \)
    \State Move \( v^* \) to the other subgraph
    \State Update \( G_1, G_2 \), and \( f(G, G_1, G_2) \)
\Until \( f(G, G_1, G_2) \) converges

\State \textbf{return} \( G_1, G_2 \)
\end{algorithmic}
\end{algorithm}

Algorithm~\ref{alg:heuristicMaxCut} provides an efficient approach to approximating the optimal solution for the constrained connected maximum cut problem with various constraints.
The time complexity of each iteration is \( O(|V| + |E|) \), considering both vertex traversal and constraint checks. With \( k \) iterations to convergence, the overall time complexity is \( O(k \times (|V| + |E|)) \). The convergence behavior is shown in Supplementary Materials, where the objective function stabilizes after several iterations.

\section{Minimal Surface Optimization}
\label{sec:MSOptimization}

In this section, we describe how the minimal surface is modeled as the classification boundary separating the two fluids and how a neural network maps spatial points to binary flow types.

\subsection{Neural Implicit MS for Space Classification} 

We construct a scalar field $f: \mathbb{R}^3 \rightarrow \mathbb{R}$, where the sign of $f(\mathbf{x})$ determines the classification of any point $\mathbf{x}$ as belonging to one of the two regions. A positive value of $f(\mathbf{x})$ indicates that $\mathbf{x}$ belongs to Region A, while a negative value indicates Region B. The zero-level set $\Omega_0=\{\mathbf{x} \in \mathbb{R}^3 \mid f(\mathbf{x})=0\}$ forms the decision boundary separating the regions, which we define as the minimal surface.

\begin{figure}[h]
    \centering
    \includegraphics[width=1\linewidth]{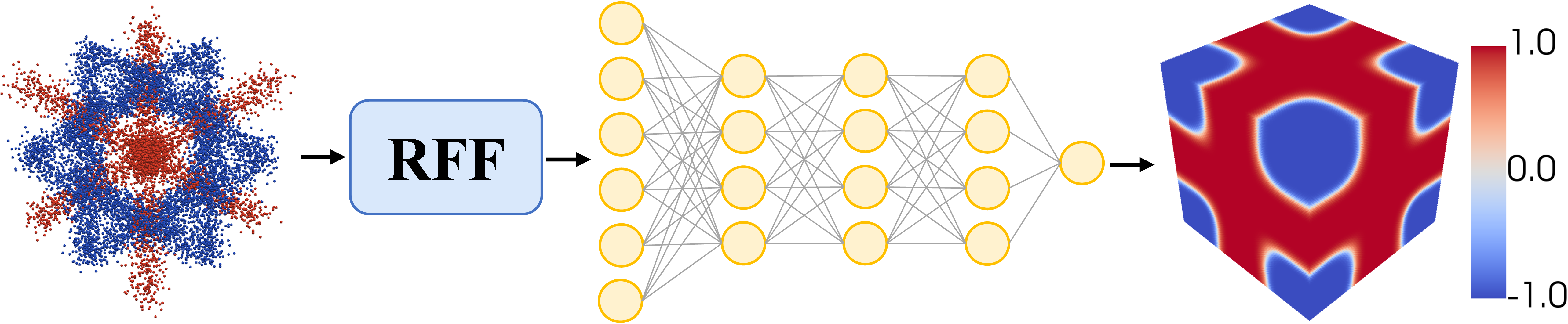}
    \caption{An overview of our network architecture. These dual skeleton sampling points are encoded by random Fourier feature (RFF) mapping, capturing high-frequency spatial information. The processed features are passed through a multi-layer perceptron (MLP) to decode the scalar field, where the sign of the output determines the spatial classification of points.}
    \label{fig:nn-structure}
\end{figure}

To approximate the scalar field $f(\mathbf{x})$, we design a neural network that combines high-dimensional feature mapping with a multi-layer perceptron (MLP). Specifically, each 3D coordinate $(x, y, z)$ is mapped into a 2048-dimensional feature vector using random Fourier mapping, which captures high-frequency spatial information. This high-dimensional representation is then processed by an MLP with three hidden layers, each containing 256 neurons and employing the softplus activation function. The final output of the network is a real-valued scalar corresponding to the input location, indicating which region the point belongs to. This pipeline is illustrated in Figure~\ref{fig:nn-structure}. This architecture effectively models both global geometric features and local smooth transitions.  % 加上网络架构引用

\begin{figure}
    \centering
    \includegraphics[width=\linewidth]{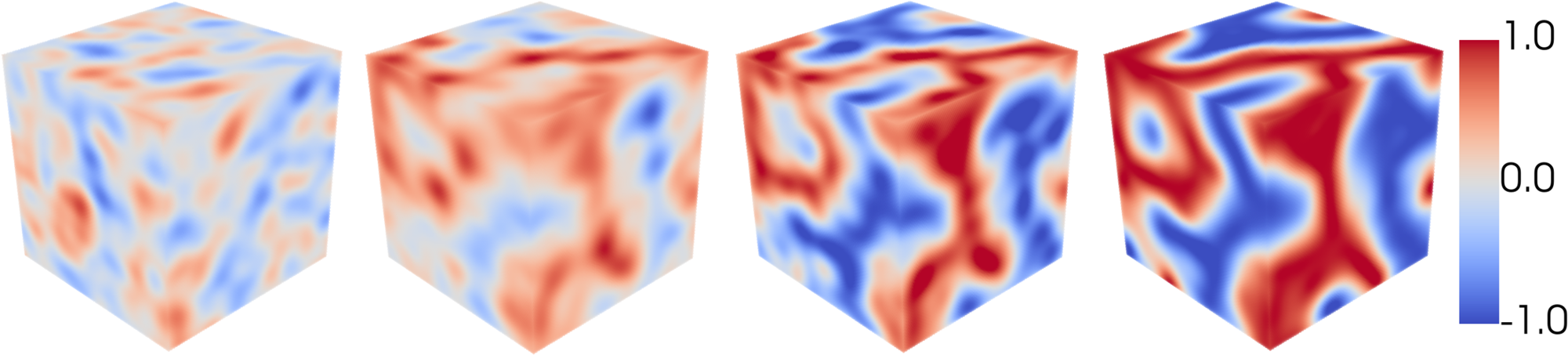}
    \leftline{ \footnotesize  \hspace{0.08\linewidth}
            (a)  \hspace{0.19\linewidth}
            (b)  \hspace{0.19\linewidth}
            (c) \hspace{0.19\linewidth}
            (d)}
    \caption{The evolution of spatial classification with increasing training iterations.}
    \label{fig:training}
\end{figure}

The network's training data are derived from the dual skeleton. Points sampled from one skeleton (denoted as $A$) are assigned a target value of $+1$, while points from the other skeleton (denoted as $B$) are assigned $-1$. These sampled points serve as direct supervision signals for the network, enabling it to learn the classification task. By encoding the skeleton information into the network's training process, the network learns the topological and positional information of the skeletons. The resulting minimal surface, serving as the classification boundary, is further refined by the smoothness regularization, ensuring a geometrically smooth and physically consistent separation of the two regions. The evolution of the scalar field \( f(\mathbf{x}) \) during training is shown in Figure~\ref{fig:training}, illustrating how the decision boundary (minimal surface) becomes progressively smoother and more aligned with the dual skeleton. 

To ensure accurate classification and generate a smooth minimal surface, we employ a combination of loss functions:

\textbf{Skeleton Loss:} This loss enforces fidelity to the skeleton labels by penalizing deviations from the target values at the sampled points:
    \begin{equation}
        \mathcal{L}_{\text{skeleton}} = \frac{1}{|A|} \sum_{\mathbf{x}_i \in A} \left| f(\mathbf{x}_i) - 1 \right| + \frac{1}{|B|} \sum_{\mathbf{x}_i \in B} \left| f(\mathbf{x}_i) + 1 \right|.
    \end{equation}

\textbf{Smoothness Loss:} To ensure that the resulting surface is smooth and area minimizing, we introduce a regularization term that penalizes the total variation of $f(\mathbf{x})$. Specifically, the smoothness loss is defined as:
    \begin{equation}
        \mathcal{L}_{\text{smooth}} = \frac{1}{|\Omega|} \sum_{\mathbf{x}_j \in \Omega} \left\| \nabla f(\mathbf{x}_j) \right\|_2,
    \end{equation}
    where $\Omega$ represents a set of uniformly sampled points within the domain.
    \wz{This formulation is theoretically grounded in the coarea formula from geometric measure theory (GMT), which establishes that minimizing the total variation of a function leads to minimizing the surface area of its level sets. In our case, minimizing $\mathcal{L}_{\text{smooth}}$ drives the zero-level set $\{ \mathbf{x} \mid f(\mathbf{x}) = 0 \}$ toward a surface with minimal area, under the topological constraints imposed by the skeleton. While this principle has been utilized in domains such as image segmentation \cite{TV17}, to the best of our knowledge, it has not been previously applied to minimal surface generation in geometric modeling.}

The total loss function combines these two components:
\begin{equation}
    \mathcal{L}_{\text{total}} = \lambda_1 \mathcal{L}_{\text{skeleton}} + \lambda_2 \mathcal{L}_{\text{smooth}},
    \label{eq:totalloss}
\end{equation}
where $\lambda_1,\lambda_2$ is a hyperparameter balancing the skeleton fidelity and the surface smoothness, empirically set $\lambda_1=5000,\lambda_2=1$ in our experiments.

By combining high-dimensional feature mapping, neural network modeling, and carefully designed loss functions, this framework enables the generation of smooth minimal surfaces tailored to the geometry of dual skeletons. The learned minimal surface not only achieves the separation of the two regions but also ensures optimal smoothness, aligning with the design goals of reducing flow resistance and maximizing thermal efficiency.

\subsection{Gaussian Noise-Based Data Augmentation}

Training subject to Eq.~\ref{eq:totalloss} alone may produce degenerate surfaces that strictly adhere to the skeleton points, which corresponds to an overfitting classifier.
To mitigate overfitting, we follow standard machine learning practice and introduce Gaussian perturbations to the points sampled from each skeleton. This noise expands the region around each skeleton, allowing the classifier to generalize better by considering local neighborhoods as belonging to the same class. As a result, the network perceives a "tube-like" region around the skeletons, preventing the decision boundary from degeneracy.

\begin{figure}[h]
    \centering
    \includegraphics[width=1\linewidth]{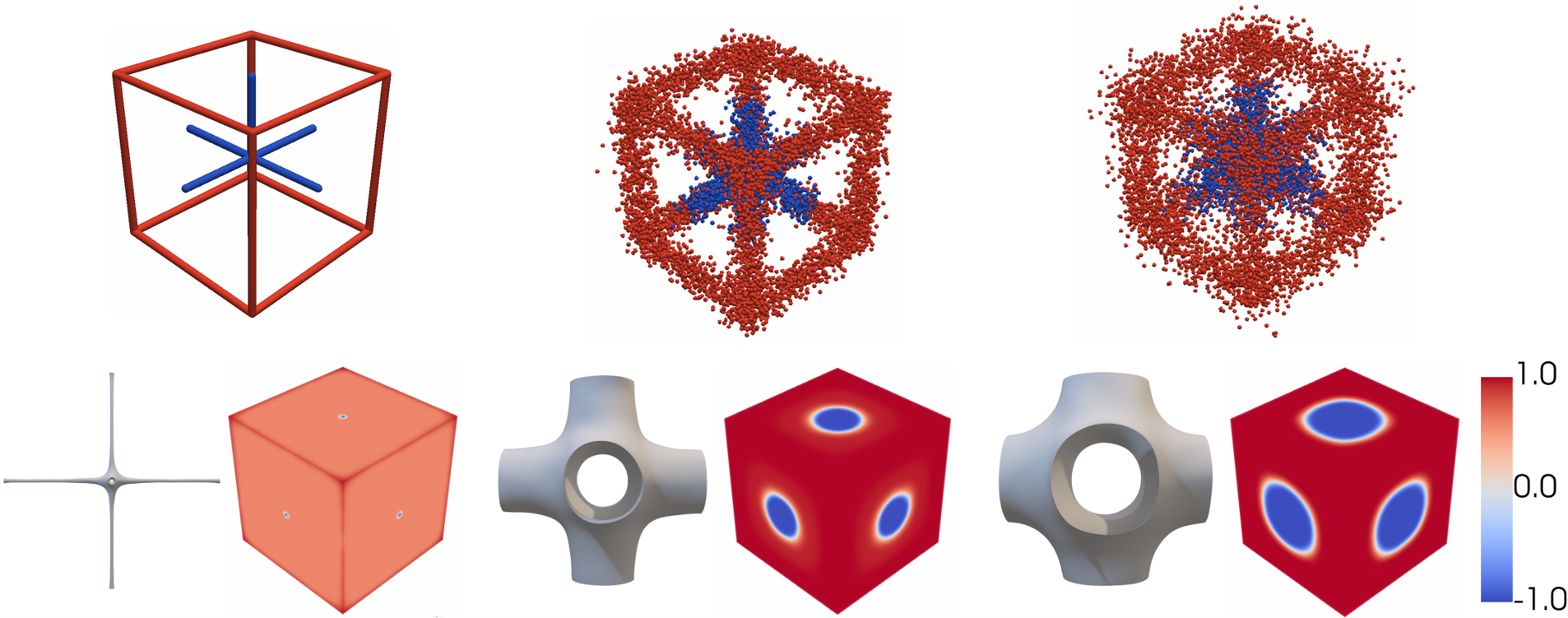}
    \leftline{ \footnotesize  \hspace{0.1\linewidth}
        (a) \hspace{0.28\linewidth}
        (b) \hspace{0.28\linewidth}
        (c)}
    \caption{Data augmentation with Gaussian noise. (a) Overfitting without noise, resulting in a decision boundary tightly clinging to the skeleton points. (b) Gaussian noise with standard deviation $\sigma=0.05$, creating a narrow “tube-like” region. (c) Gaussian noise with standard deviation $\sigma=0.08$, significantly broadening the neighborhood.}
    \label{fig:Gaussian}
\end{figure}

Given a sampled point \( \mathbf{x} = (x, y, z) \) on the skeleton, we generate a perturbed point:
\begin{equation}
    \mathbf{x^\prime} = \mathbf{x} + \boldsymbol{\eta},
\end{equation}
where \( \boldsymbol{\eta} \sim \mathcal{N}(\mathbf{0}, \sigma^2 \mathbf{I}) \). The variance \( \sigma^2 \) controls the extent of the perturbation, balancing data coverage and topology preservation. 
Figure~\ref{fig:Gaussian} shows how this noise broadens the distribution around the skeletons, helping to avoid overfitting.

%% file: 5-results.tex
\section{Results and Discussion}
\label{sec:results}

We evaluate our method by presenting the results of both minimal surface properties and heat transfer performance. 

\subsection{Minimal Surface Properties}

In this section, we validate the minimal surface properties of our method. 
All of our models are trained using the Adam optimizer on a single NVIDIA GeForce RTX 3090 GPU. 
We complete the constrained connected maximum cut algorithm in seconds. 
We train each model for $5.12\times10^4$ iterations (about 70 minutes) with a learning rate of $3\times10^{-5}$, sampling $32,768$ points from the freeform domain. 
\wz{We regularize the design domain within a unit cubic space and apply Gaussian noise with a standard deviation $\sigma=0.002$ at each step.  
}
The surface meshes are extracted by Marching Cubes from the neural implicit field at a resolution of $256^3$. 
\wz{A breakdown of computational time performance is provided in the Supplementary Materials.}

We demonstrate the key characteristics of minimal surfaces, including low mean curvature \cite{CohenSteinerMorvan03} and high geometric regularity. These properties are essential for ensuring that the surfaces are not only geometrically accurate but also suitable for practical applications, such as enhancing heat transfer efficiency in heat exchangers.

\paragraph{Ablation Study on Smoothness Loss} 

To validate the effect of the smoothness term in the generation of minimal surfaces, we perform an ablation study by removing the smoothness term from the loss function. 
This comparison evaluates the differences in surface morphology and the scalar field used for spatial classification between surfaces generated with and without the smoothness loss.

\begin{figure}[t]
        \centering
        \includegraphics[width=1\linewidth]{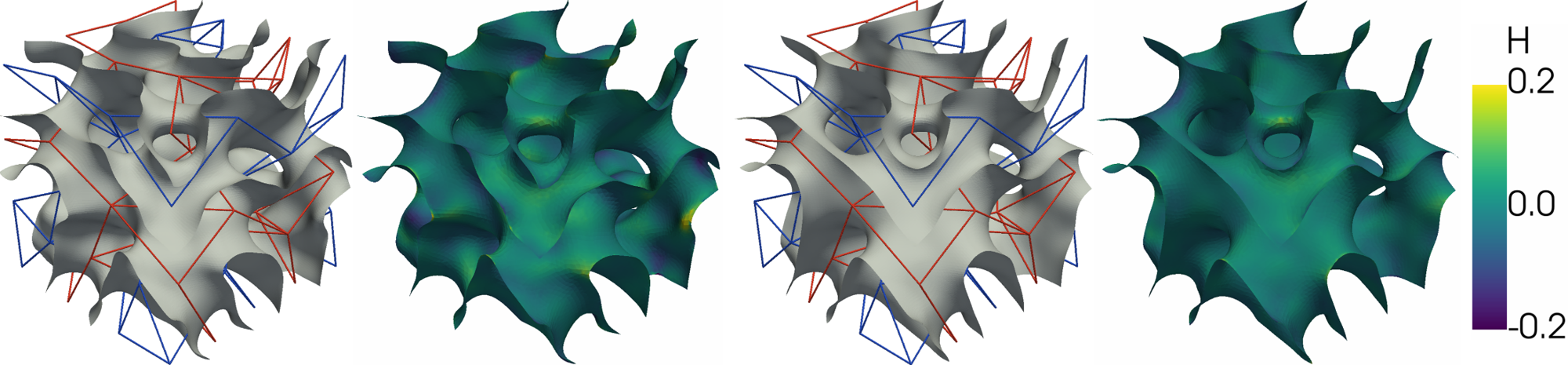}
        \leftline{ \footnotesize  \hspace{0.08\linewidth}
        (a) without smoothness loss \hspace{0.13\linewidth}
        (b) with smoothness loss }
        \caption{Ablation study on the smoothing term loss for minimal surface generation. Result with the smoothness loss. The mean curvature of (b) is more uniform and closer to 0.}
        \label{fig:smoothness-ablation}
\end{figure}

Figure~\ref{fig:smoothness-ablation} demonstrates that removing the smoothness loss leads to visibly rougher structures with more irregular curvature distribution. 
In Figure~\ref{fig:smoothness-ablation}(a), without the smoothness loss, the structure exhibits sharp features and unnecessary twists. In contrast, Figure~\ref{fig:smoothness-ablation}(b), with the smoothness loss, shows a smoother and more continuous structure. This comparison highlights the crucial role of the smoothness term in ensuring the generation of high-quality minimal surfaces.

\begin{table}[t]
    \centering
        \caption{Comparison of geometric properties, including mean curvature (H)  and surface area, between the minimal surface (S) generated by our method and the SDF-based equidistant surface after Laplacian smoothing.}
    \scalebox{1}{
    \begin{tabular}{>{\centering\arraybackslash}m{1.0cm}>{\centering\arraybackslash}m{3.5cm}>{\centering\arraybackslash}m{2.8cm}}
        \toprule%第一道横线
         Method & SDF + Laplacian smoothing & Ours  \\
        \midrule%第二道横线
         \text{S} & \includegraphics[width=0.17\textwidth]{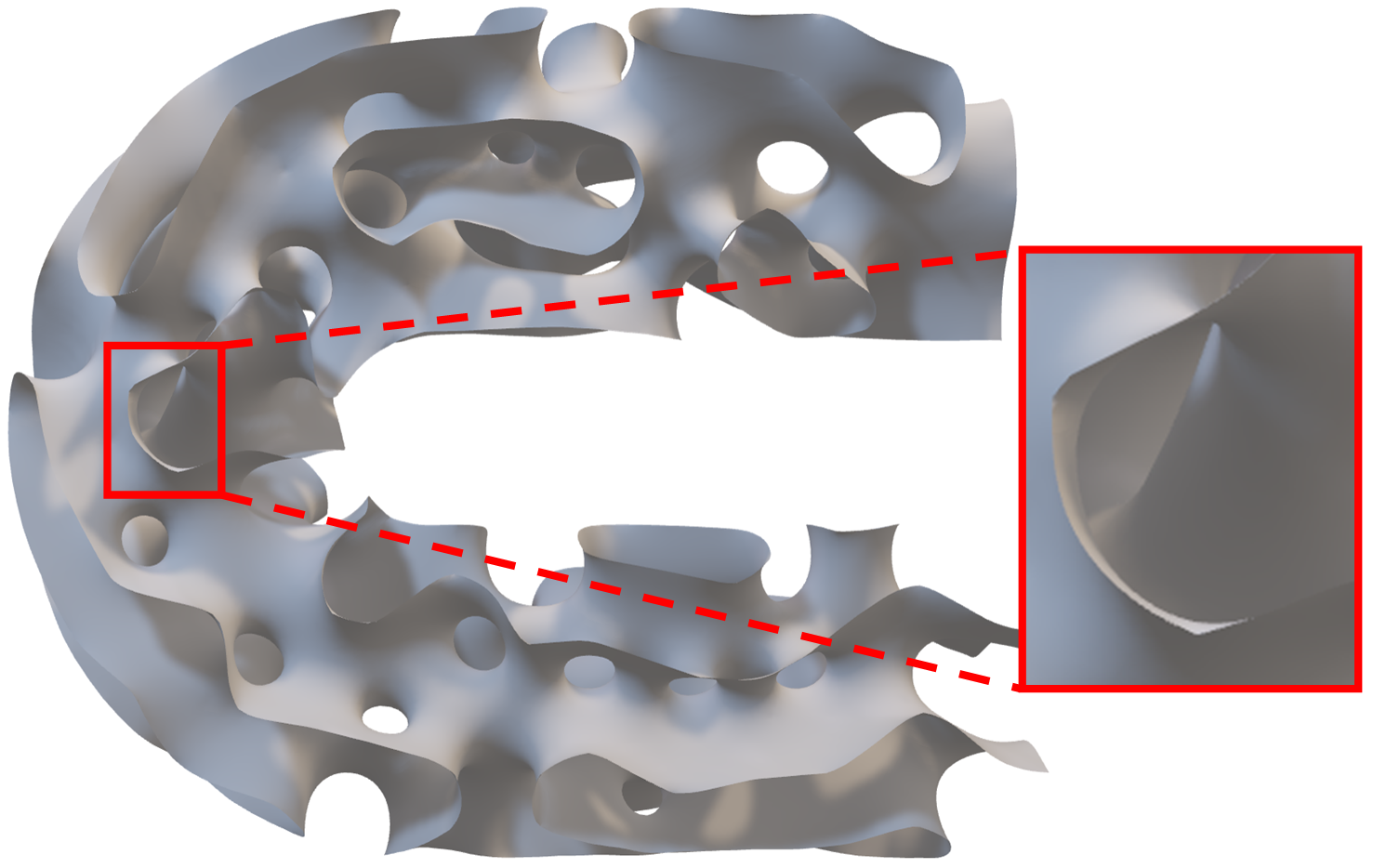} & \includegraphics[width=0.17\textwidth]{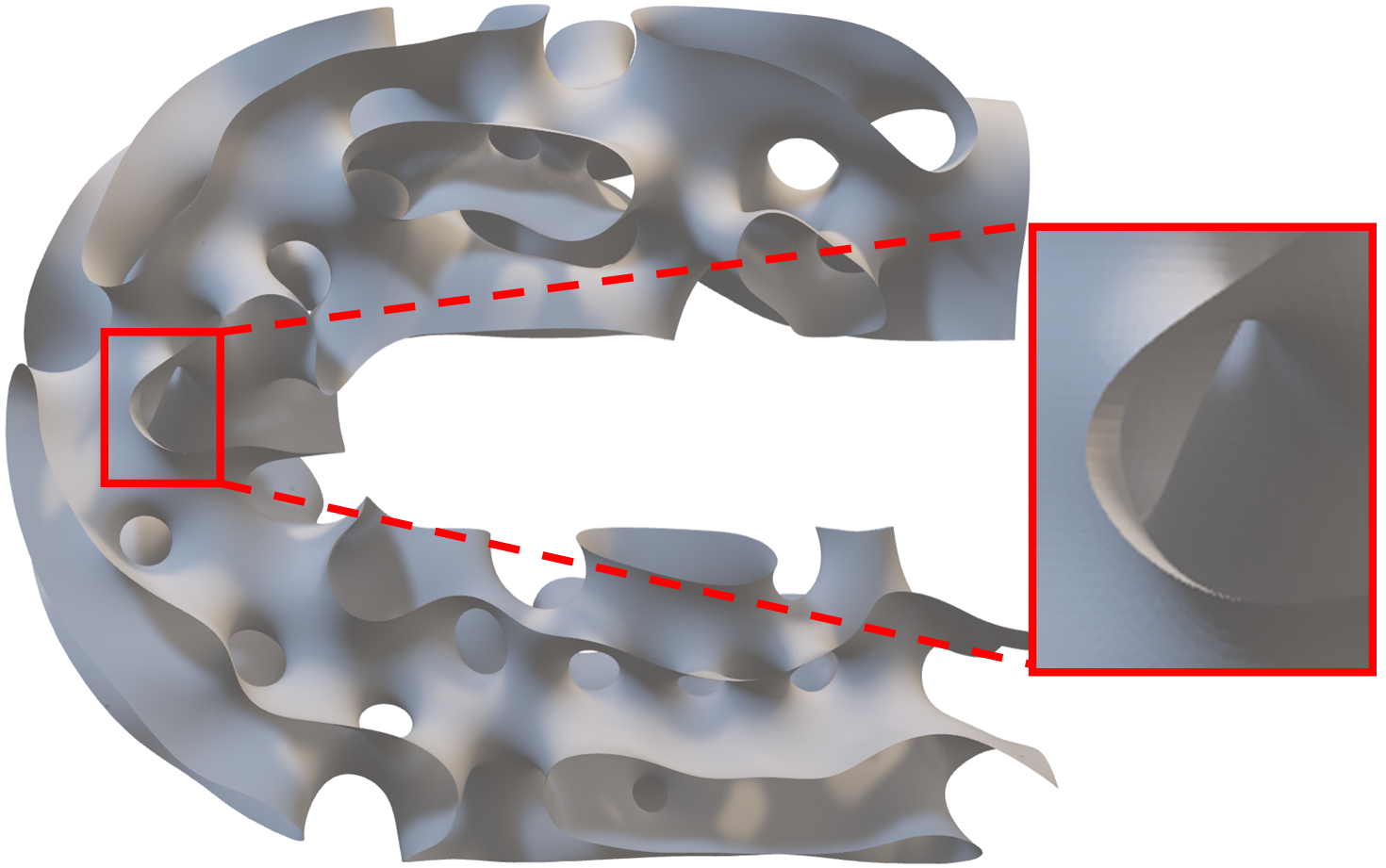} \\
         \text{H} & [-0.020, 0.269] & [-0.011, 0.074] \\
         \text{Area} & 3.51 & 3.43 \\
       \bottomrule%第三道横线
    \end{tabular}
    }
    \label{tab:compare-SDF}
\end{table}

\paragraph{Comparison with Surface Smoothing based on SDF} 

To assess the optimality of the generated minimal surfaces, we compare our results with the equidistant surface between the dual skeletons, derived from signed distance fields (SDF), followed by Laplacian smoothing. The equidistant surface has the desired topology that separates two flow channels and an improved geometry smoothed by mean curvature flow, thus serving as a good baseline for surface generation.
Since the SDF-based method and our approach share the same dual skeleton, we evaluate the resulting minimal surfaces in terms of surface area and mean curvature, as shown in Table~\ref{tab:compare-SDF}.
The table shows that the SDF-based surface has sharp geometric features, particularly in the red-box region, while our method generates a smoother surface. 
This visual observation is further supported by the curvature range and surface area, which are smaller and more optimized for the minimal surface produced by our approach than those of the smoothed SDF-based surface.

\paragraph{Comparison with Classical Minimal Surfaces} 

To evaluate the ability of our method to generate minimal surfaces, we test it by reconstructing the classical minimal surfaces such as Schwarz P and Gyroid. These surfaces are often used as benchmarks in minimal surface research and have been extensively studied for their geometrical properties.

We begin by designing skeletons that closely match the topology of these classical minimal surfaces. 
Using these skeletons as input, we generate minimal surfaces by our method and compare the surface area and curvature distribution with the classical minimal surfaces. 
As shown in Table~\ref{tab:compare-TPMS}, the generated surfaces exhibit similar curvature distributions and surface areas compared to the original minimal surfaces, validating the effectiveness of our method in accurately reconstructing classical minimal surfaces.

\begin{table}[htbp]
    \centering
    \setlength{\tabcolsep}{1pt}
    \caption{Comparison of classical minimal surfaces  with similar structures generated by our method. The table highlights differences in curvature distributions (\(H\in[-0.1,\ 0.1]\)), surface area, and mesh properties (number of vertices \#V and faces \#F).}
 %   \resizebox{0.9\linewidth}{!}{ % 调整表格宽度为90%栏宽
    \begin{tabular}{>{\centering\arraybackslash}m{0.7cm}
    >{\centering\arraybackslash}m{1.2cm}
    >{\centering\arraybackslash}m{1.2cm}
    >{\centering\arraybackslash}m{1.2cm}
    >{\centering\arraybackslash}m{1.2cm}
    >{\centering\arraybackslash}m{1.2cm}
    >{\centering\arraybackslash}m{1.2cm}
    >{\centering\arraybackslash}m{0.2cm}} % 第一列设为可换行并居中
        \toprule % 第一道横线
          & {\small P} & {\small P-Ours} & {\small IWP} & {\small IWP-Ours} & {\small G} & {\small G-Ours} \\ % 添加Bar列
        \midrule % 第二道横线
        \multirow{2}{*}{\makecell{H}}
        % \multirow{2}{*}{Mean \\ curvature}
        & \includegraphics[width=1\linewidth]{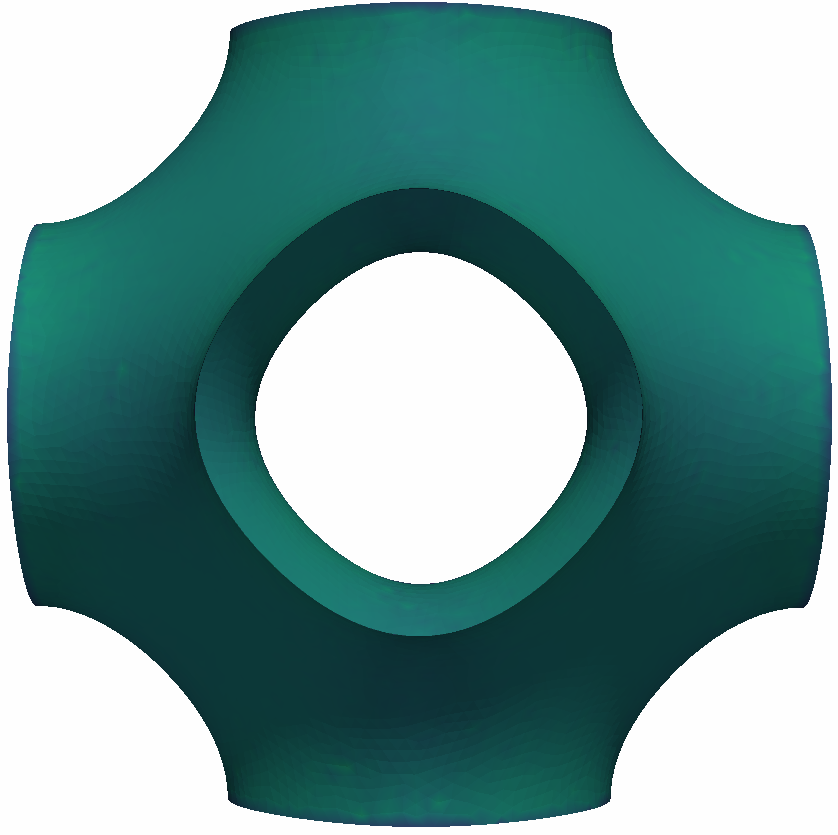} 
        & \includegraphics[width=1\linewidth]{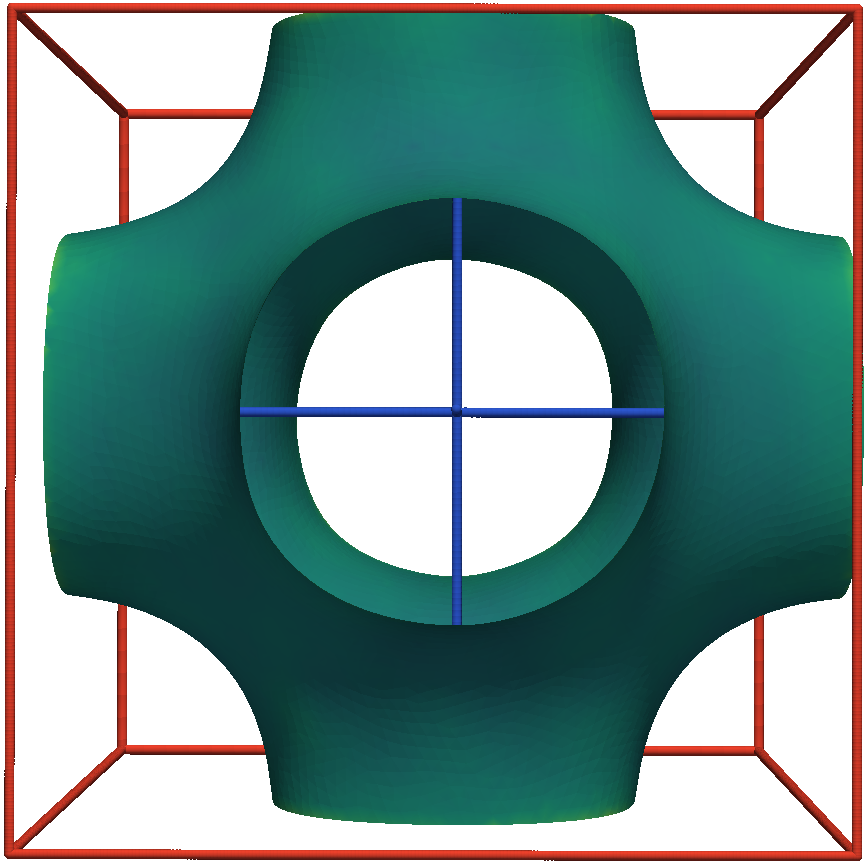} 
        & \includegraphics[width=1\linewidth]{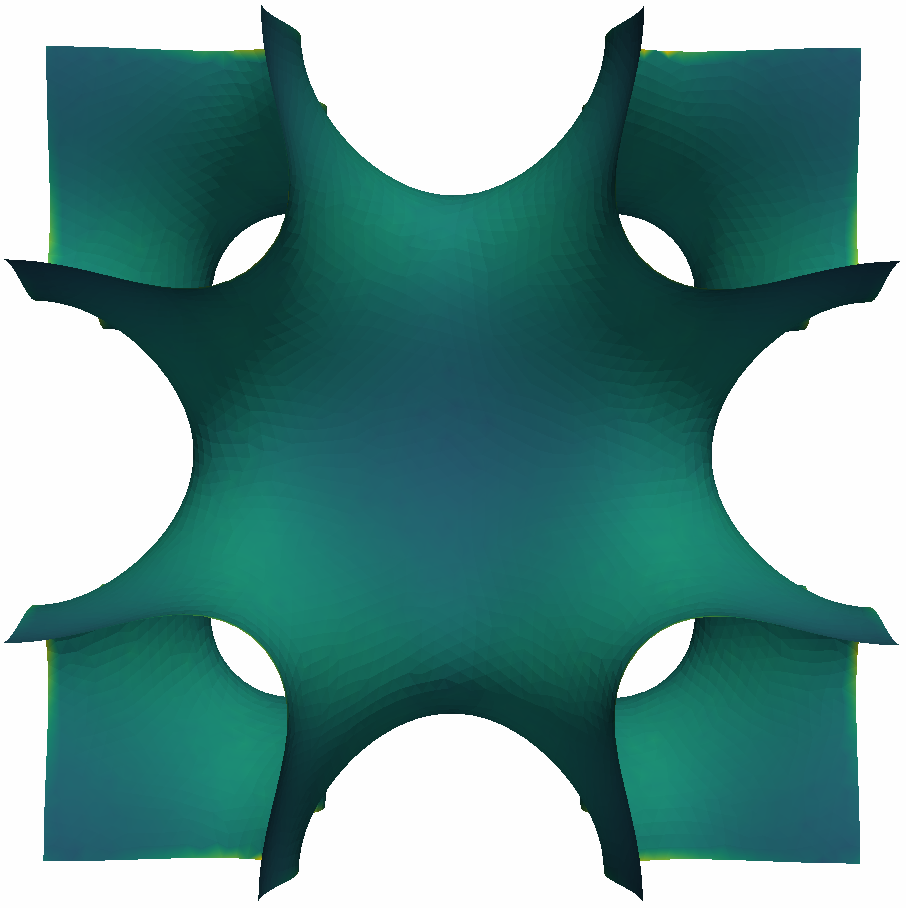} 
        & \includegraphics[width=1\linewidth]{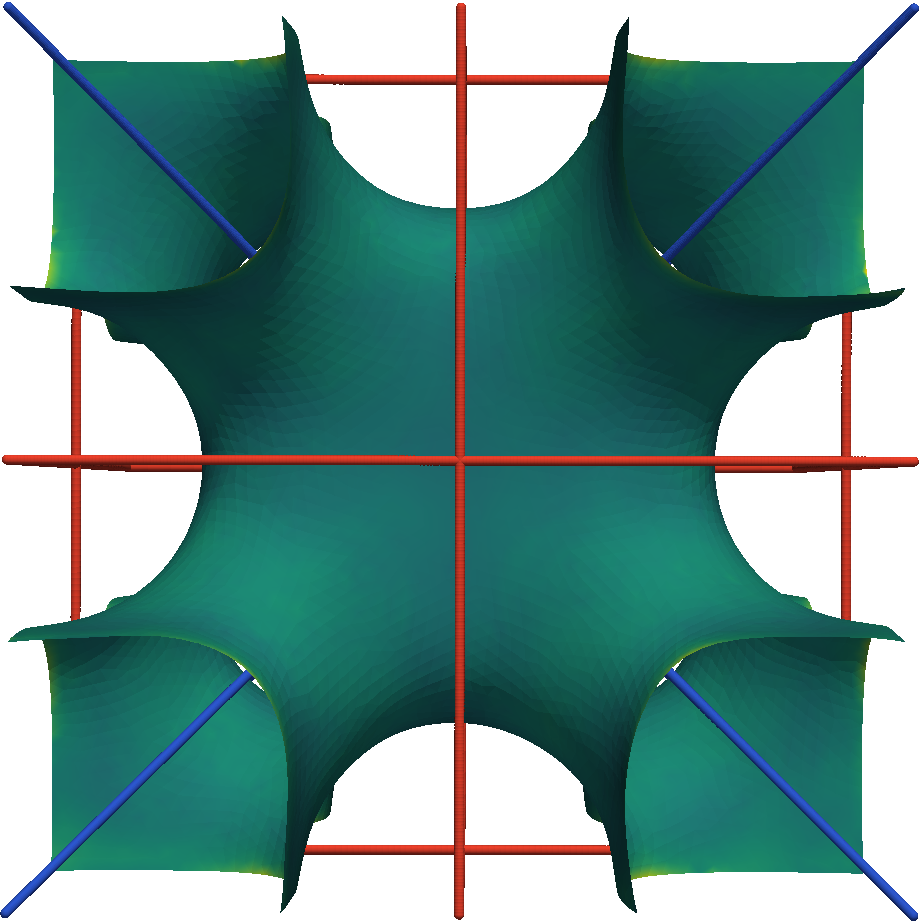} 
        & \includegraphics[width=1\linewidth]{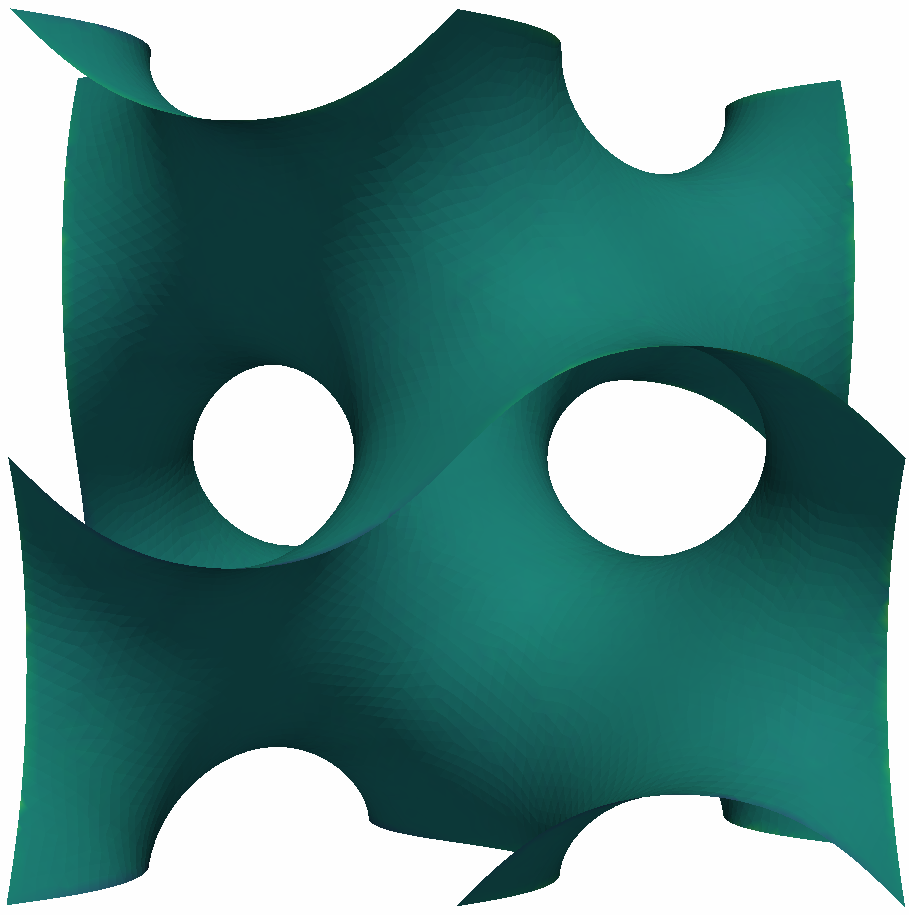} 
        & \includegraphics[width=1\linewidth]{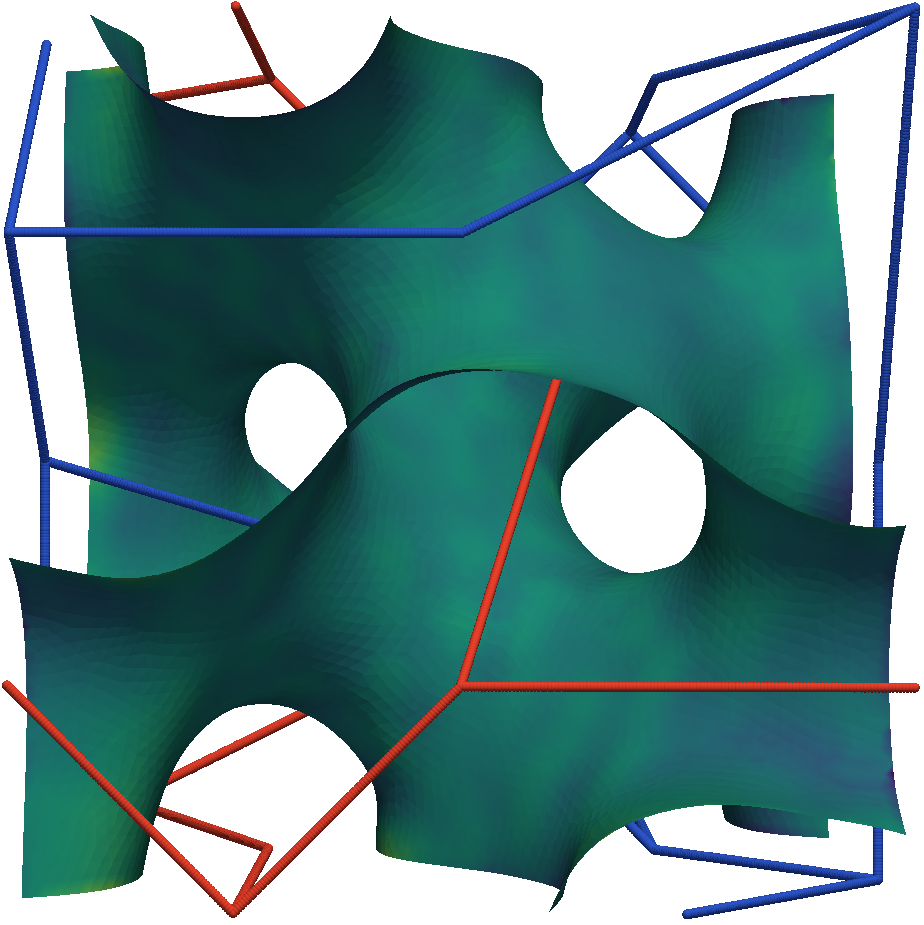} 
        & \includegraphics[width=2\linewidth]{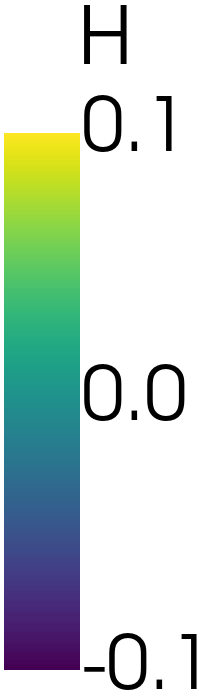} \\ 
        
        & \includegraphics[width=1\linewidth]{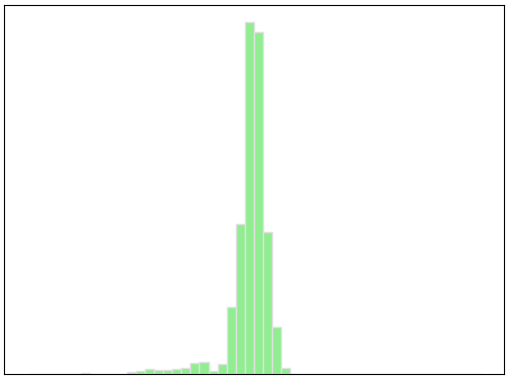} 
        & \includegraphics[width=1\linewidth]{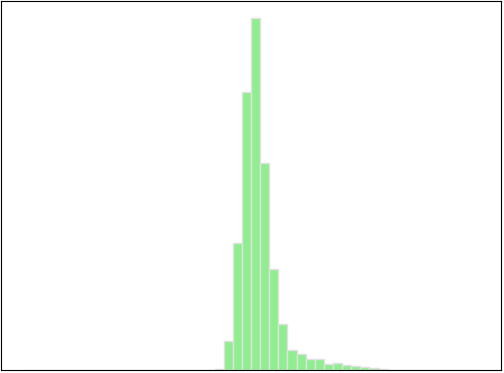} 
        & \includegraphics[width=1\linewidth]{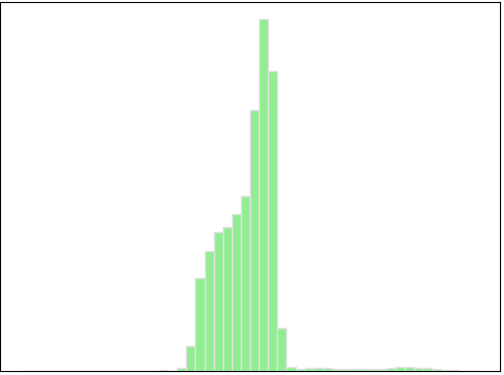} 
        & \includegraphics[width=1\linewidth]{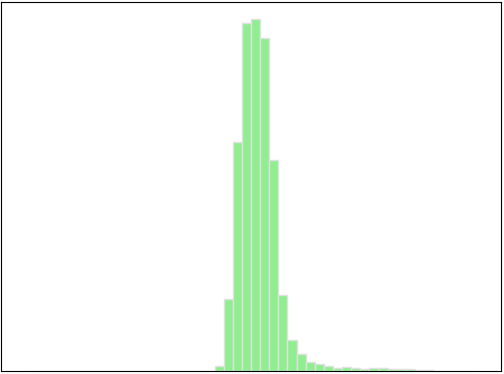} 
        & \includegraphics[width=1\linewidth]{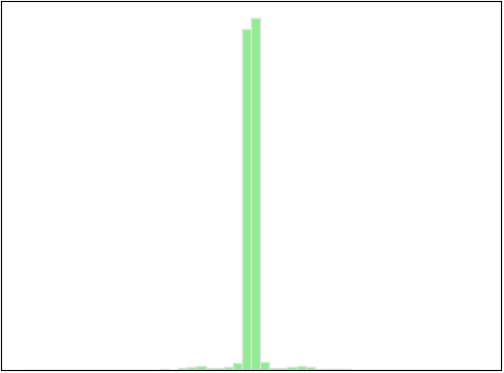} 
        & \includegraphics[width=1\linewidth]{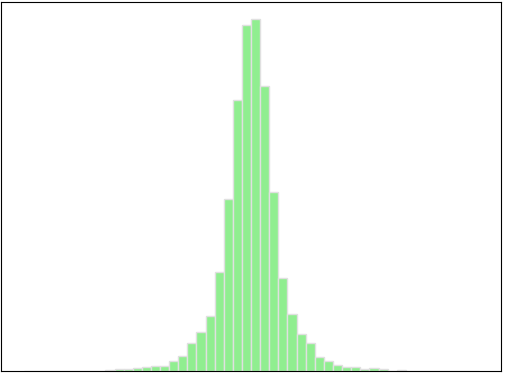} 
        & \\ % 为空，保持行对齐
        
        {Area}  & \centering 3.79 & \centering 3.79 & \centering 5.73 & \centering 5.63 & \centering 4.99 & \centering 4.72 & \\
        {\#V} & 8679 & 8707 & 13127 & 12843 & 11512 & 10865 & \\ % 最后一列留空
        {\#F} & 16870 & 16922 & 25476 & 24941 & 22249 & 21068 & \\ % 最后一列留空
       \bottomrule % 第三道横线
    \end{tabular}
    %}
    \label{tab:compare-TPMS}
\end{table}

\subsection{Heat Exchange Performance}

In this section, we validate the heat exchange performance of our method by comparing them with the mostly used Gyroid and topologically optimized geometries by CFD simulation. 
The evaluation consists of two key metrics: pressure drop and temperature difference, which directly reflect the efficiency of heat transfer in practical applications.
The CFD simulations are carried out on a computer with an AMD Ryzen Threadripper PRO 5995WX 64-Cores CPU, 512 GB RAM. All our simulation experiments converge within 600 iterations using ANSYS 2023 R1.

\paragraph{Comparison with Uniform Gyroid Structures}

\begin{figure*}[htpb]
    \centering
    \includegraphics[width=\linewidth]{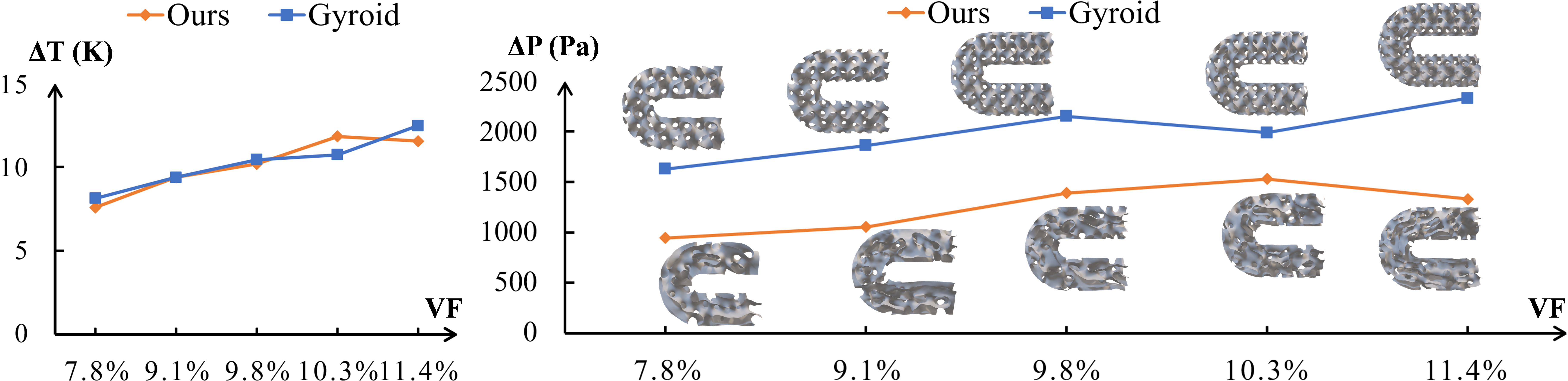}
        \leftline{ \footnotesize  \hspace{0.15\linewidth}
            (a)  \hspace{0.48\linewidth}
            (b) }
    \caption{Comparison of the temperature difference ($\Delta T$) (a) and pressure drop ($\Delta P$) (b) between Gyroid and our structure in different volume fractions (VF). }
    \label{fig:PT-volume}
\end{figure*}

\begin{figure*}[htpb]
    \centering
    \includegraphics[width=1\linewidth]{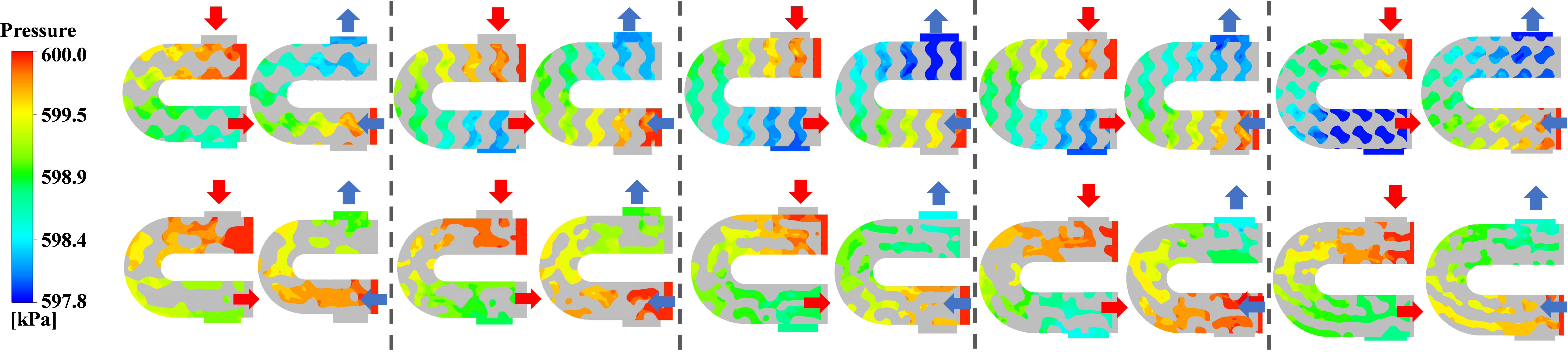}
    \caption{Comparison of the pressure drop between Gyroid (upper row) and our structure (lower row) in different volume fractions. The volume fractions are $7.8\%,\ 9.1\%,\ 9.8\%,\ 10.3\%,\ 11.4\%$, respectively, from left to right}
    \label{fig:pressure-drop}
\end{figure*}

The Gyroid structure is widely used for its high surface area and efficient fluid distribution properties, making it a suitable benchmark for heat exchanger design. 
% \wz{Gyroid has generally shown to outperform other commonly used TPMS designs in thermal performance\cite{Kaur2021,Gado2024}.} 
\wz{We use the Gyroid structure as the baseline because it is the most widely adopted TPMS in heat exchanger applications and has been shown to outperform other commonly used TPMS designs such as Primitive and Diamond in terms of thermal performance and pressure drop balance~\cite{Kaur2021,Gado2024}.
While conventional straight tube or core–shell configurations are feasible baselines, their inherently different structural characteristics—such as lower interface area and lack of integrated dual-channel flow—make them less suitable for direct comparison under our freeform and dual-fluid optimization setting.}
We simulate heat exchange in both the U-shape structure generated by our approach with varying numbers of sample points, specifically 300, 400, 500, 600, and 800 and the Gyroid-filled U-shape under identical operating conditions, with both structures having the same volume fraction. 
The comparison is presented in terms of pressure drop ($\Delta P$) and temperature difference ($\Delta T$).

Figure~\ref{fig:PT-volume} shows the experimental results for $\Delta T$ and $\Delta P$ with varying volume fractions. 
In Figure~\ref{fig:PT-volume}(a), our method maintains a temperature difference close to the Gyroid, demonstrating that the heat transfer efficiency is preserved. Meanwhile, Figure~\ref{fig:PT-volume}(b) indicates that our structure exhibits a lower pressure drop compared to the Gyroid, suggesting a more efficient fluid flow. Furthermore, Figure ~\ref{fig:pressure-drop} provides a detailed comparison of the pressure drop values for these experimental groups, highlighting the lower flow resistance of our method.

\begin{figure}[t]
    \centering
    \includegraphics[width=\linewidth]{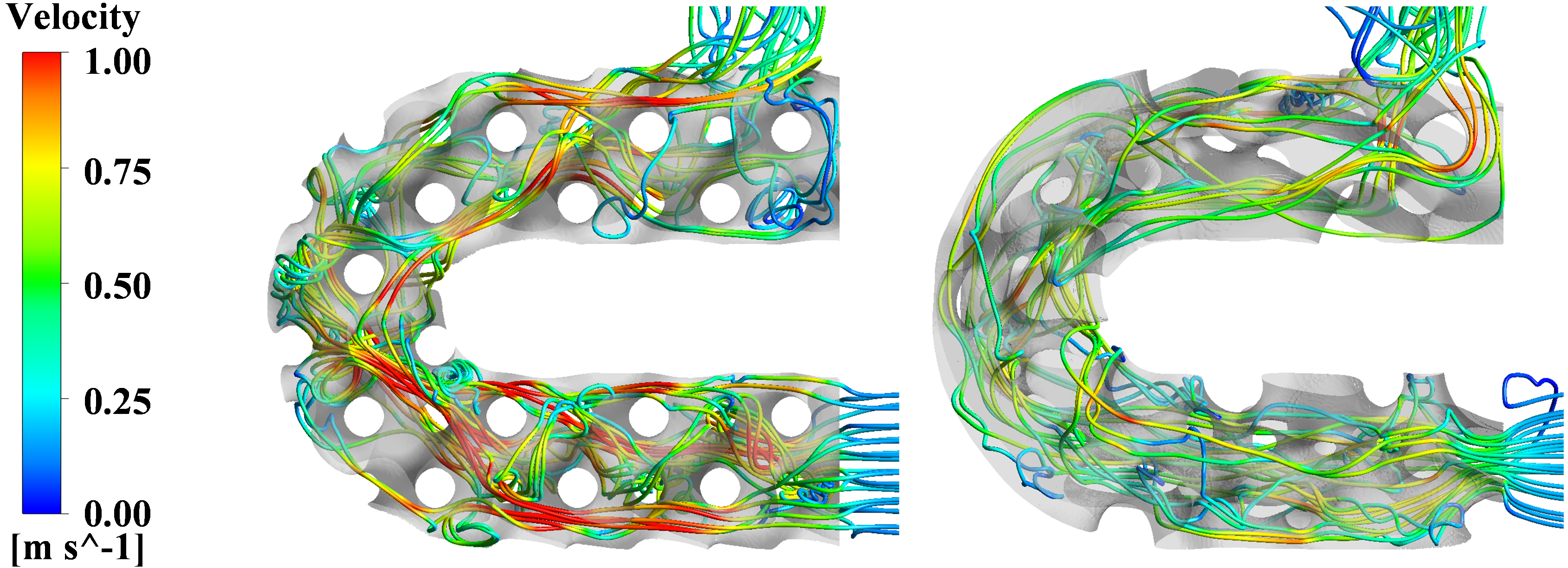}
    \caption{Comparison of streamlines between the uniform Gyroid (left) and our structure (right).}
    \label{fig:streamline}
\end{figure}

Our method achieves comparable heat transfer while substantially lowering the pressure drop. Figure~\ref{fig:streamline} shows the comparisons of the streamlines between the uniform Gyroid and our structure, emphasizing the advantage of our structure in reducing flow resistance. In the uniform Gyroid structure, obvious kinks are noticed at the U-shape corners, corresponding to large flow resistance and pressure drop. This is eliminated with our method.

\paragraph{Comparison with Optimized Gyroid Structures}

\begin{figure}[t]
    \centering
    \includegraphics[width=0.8\linewidth]{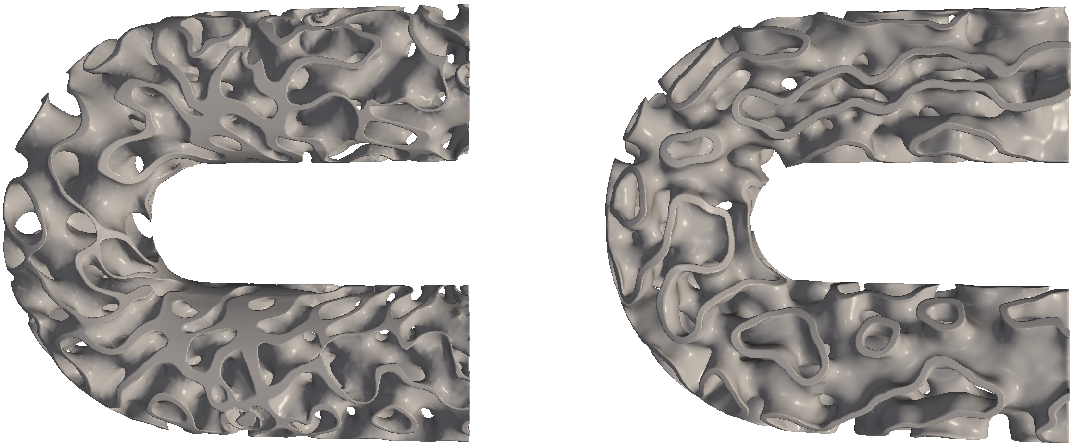}
    \caption{The geometry of topology-optimized Gyroid (left) and our structure (right) with volume fraction of 29.8\%.}
    \label{fig:compare-mesh-TO-gyroid}
\end{figure}

\begin{figure}[t]
    \centering
    \includegraphics[width=\linewidth]{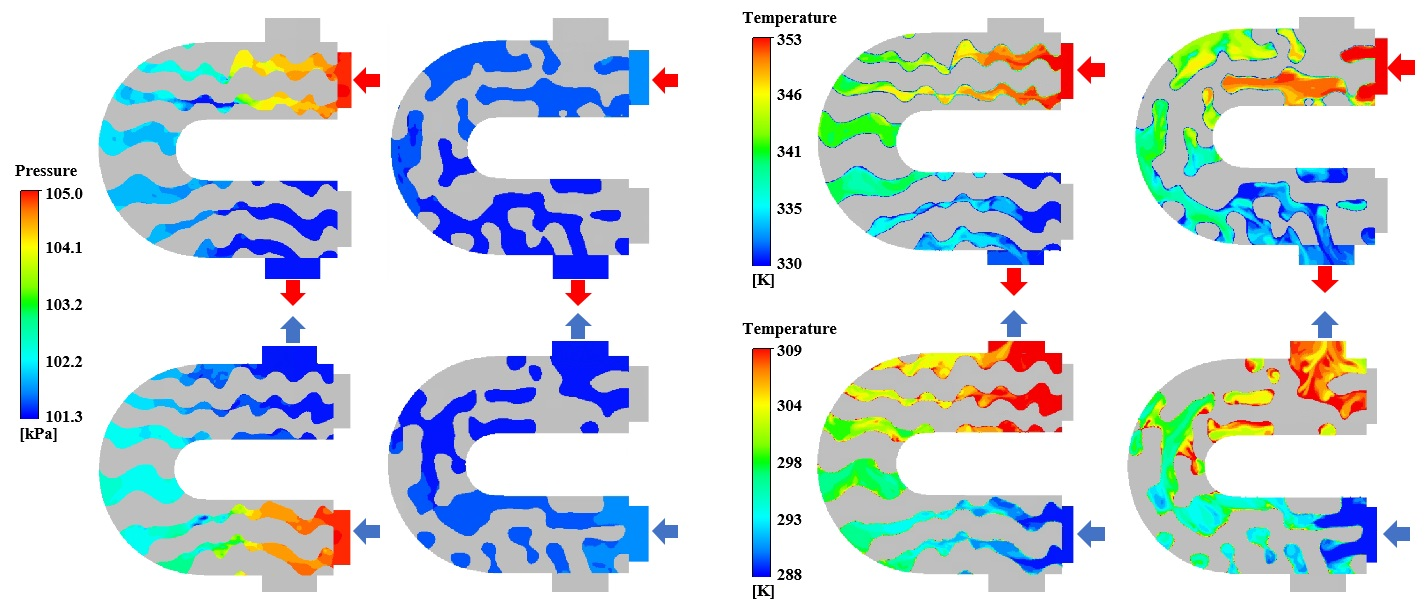}
    \leftline{ \footnotesize  \hspace{0.2\linewidth}
    (a) Pressure drop \hspace{0.3\linewidth}
    (b) Temperature }
    \caption{Comparison of heat exchange performance between the topology-optimized Gyroid (left) and our structure (right). Our method achieves approximately $92\%$ $\Delta T$ of the topology-optimized Gyroid, while the $\Delta P$ is only $29\%$.}
    \label{fig:compare-to-pt}
\end{figure}

We further compare with a topology-optimized Gyroid structure (TO Gyroid) using the same U-shape shell geometry (Figure~\ref{fig:compare-mesh-TO-gyroid}) and simulation parameters as described in \cite{Jiang2023}. 
The results, illustrated in Figure~\ref{fig:compare-to-pt}, reveal a trade-off between thermal performance and flow resistance. Although the $\Delta T$ of our method achieves approximately $92\%$ of the TO Gyroid, the $\Delta P$ is only $29\%$ of the TO Gyroid. This substantial reduction in pressure drop highlights the efficiency of our minimal surface approach in promoting smooth fluid flow.
Our method prioritizes pressure drop reduction while maintaining competitive heat transfer efficiency, our approach especially suitable in cases where energy efficiency and operational cost are critical considerations. % smaller pressure drop

\paragraph{More Freeform Heat Exchanger Results}

We evaluate the performance of our method in two heat exchangers with more freeform design space: a zigzag-shaped shell and a multi-bend pipe shell. 
These geometries require the internal structure to adapt to irregular external boundaries. 
We compare our method with the Gyroid-filled structure with identical operating conditions, ensuring the same volume fraction for both designs.

\begin{figure}[t]
    \centering
    \includegraphics[width=1\linewidth]{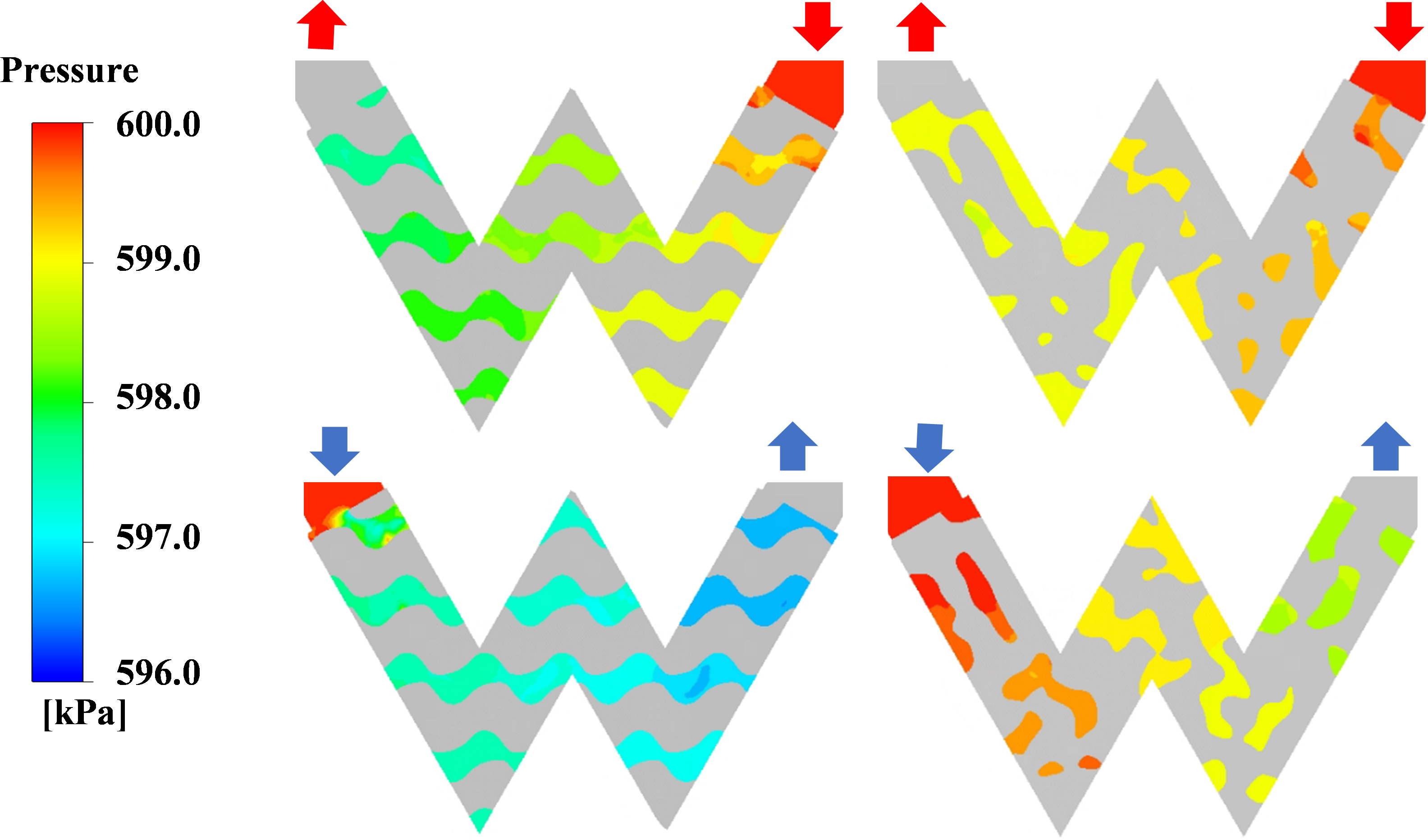}
    \caption{Cross-sectional pressure drop distributions for zigzag shape, comparing the minimal surface heat exchanger designed by the Gyroid structure (left) with our method (right). Our method provides a comparable heat transfer performance with a lower pressure drop.}
    \label{fig:freeform-comparison}
\end{figure}

\begin{figure}[t]
    \centering
    \includegraphics[width=1\linewidth]{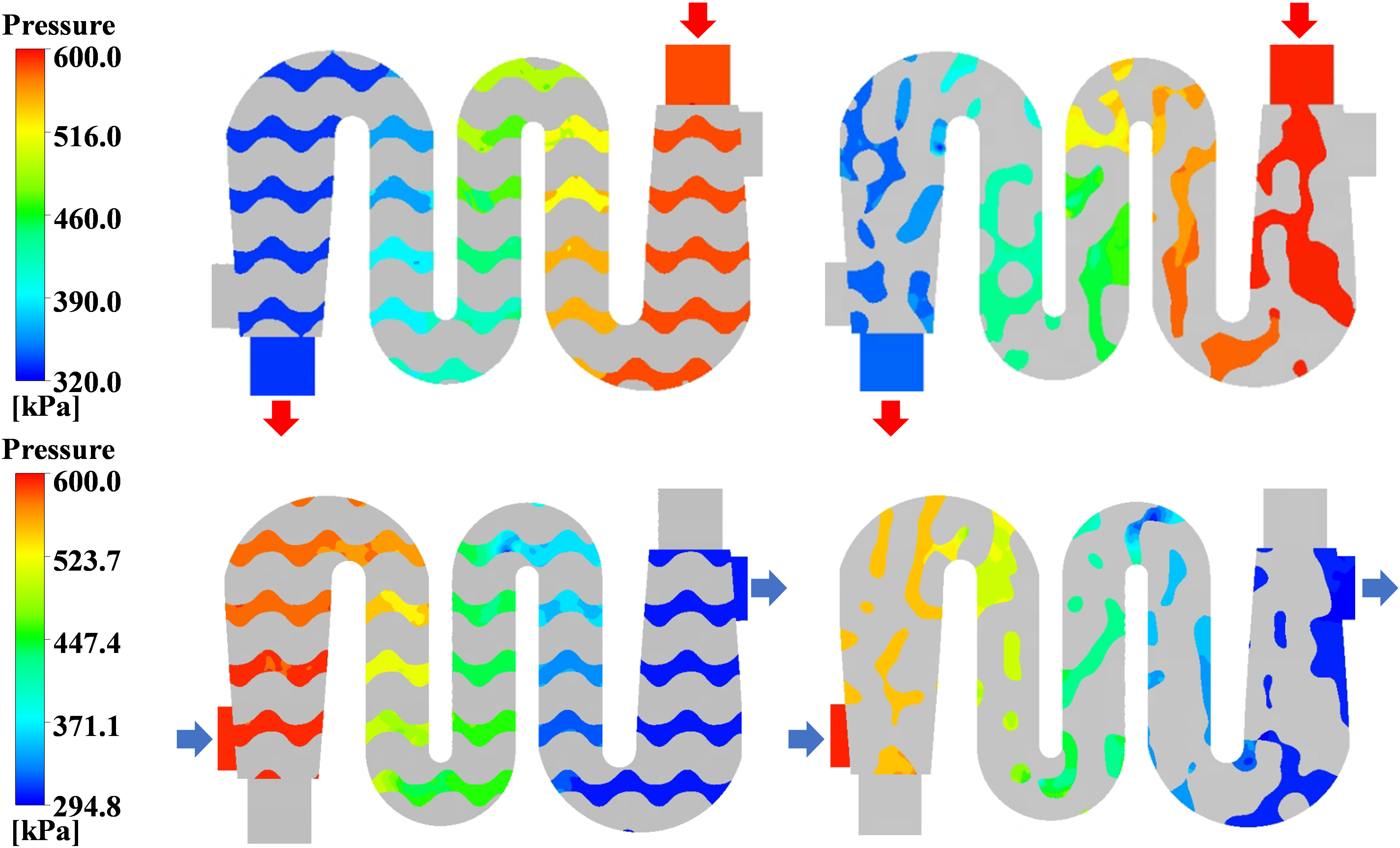}
    \caption{Cross-sectional pressure drop distributions for multi-bend shape, comparing the minimal surface heat exchanger designed by the Gyroid structure (left) with our method (right). Our method provides a comparable heat transfer performance with a lower pressure drop.}
    \label{fig:multi-bend-shell-comparison}
\end{figure}

As shown in Figures~\ref{fig:freeform-comparison} and \ref{fig:multi-bend-shell-comparison}, our method achieves a comparable temperature difference compared to the Gyroid structure both in the zigzag and multi-bend cases. At the same time, our minimal surface designs exhibit significantly reduced pressure drop. Specifically, for the zigzag case, the $\Delta T$ and $\Delta P$  are 14.36 K and 1809 Pa for our method, compared to 14.87 K and 3181 Pa for the Gyroid structure. 
Similarly, for the multi-bend case, our method achieves $\Delta T$ of 15.95 K and a $\Delta P$  of 282 kPa, while the Gyroid structure shows a $\Delta T$ of 15.55 K and a $\Delta P$ of 296 kPa. By combining high thermal performance with reduced flow resistance, our approach enables the design of efficient heat exchangers for freeform geometries.

%% file: 6-conclusion.tex
\section{Conclusion and Future work}
\label{sec:conclusion}

This paper introduces \emph{DualMS}, a framework for optimizing dual-channel minimal surfaces in heat exchangers with freeform shapes. 
Unlike traditional triply periodic minimal surfaces (TPMS), DualMS directly optimizes minimal surfaces for two-fluid heat exchangers, addressing flow constraints and minimizing pressure drops. By formulating the heat exchange maximization as a constrained connected maximum cut problem and using a neural network to automatically determine surface boundaries and classify flow skeletons, DualMS offers enhanced flexibility in surface topology. 
The results show that DualMS achieves superior thermal performance, with lower pressure drops and comparable heat exchange rates to TPMS under the same material cost, advancing heat exchanger design for freeform applications.

% future work
The proposed framework has certain limitations that will be addressed in future work. 
First, the efficiency of the minimal surface optimization network can be improved, as suggested in \cite{Mueller2022}. 
Second, self-supporting constraints \cite{Xing2024} were not incorporated into the optimization, potentially leading to overhangs that require additional support structures in certain printing processes. While advancements in laser additive manufacturing \cite{velo3d-supportfree}, multi-axis additive manufacturing~\cite{TONG2024104508}, and hybrid manufacturing \cite{Zhong2023} are reducing these fabrication constraints, investigating the integration of such constraints into the minimal surface solving remains valuable.
Looking ahead, we plan to develop an end-to-end heat exchanger design system incorporating application-specific constraints \cite{Xu2023CVM,xue2025mind}, such as flow density and volume. Our framework can also be extended to the design of three- or more fluid heat exchangers \cite{Wei2025}.

%% file: appendix.tex
\appendix

%%%%%%%%%%%%%%%%%%%%%%%%%%%%%%%%%
\section{Computational Time Breakdown}

Table~\ref{tab:Computational-Time} presents the computational cost breakdown of our framework, including graph initialization (GI), constrained connected maximum cut (CCMT), and minimal surface generation (MSG).

\begin{table}[h]
    \centering
    \caption{Computational time breakdown for different point counts.}
    \begin{tabular}{cccc}
        \toprule
        Points Number & GI [second] & CCMT [second]  & MSG [minute]\\
        \midrule
        300&3.1668&0.0121&67.0512\\
        400&3.4427&0.0180&63.4799\\
        500&2.9378&0.0326&68.5413\\
        600&3.1368&0.0432&68.9918\\
        800&3.7321&0.0647&66.2961\\
        \bottomrule
    \end{tabular}
    \label{tab:Computational-Time}
\end{table}

%%%%%%%%%%%%%%%%%%%%%%%%%%%%%%%%%
\section{Convergence of the maximum cut algorithm}

In this section, we demonstrate the convergence behavior of our algorithm and evaluate the changes in objective function values with different numbers of sampling points.

Since the objective function is based on the number of cut edges, and this number cannot exceed the total number of edges in the original graph, the objective function has an upper bound. As a result, the algorithm is guaranteed to converge, as the search for an optimal solution will not lead to an infinite increase in the objective function value.

\begin{figure}[h]
    \centering
    \includegraphics[width=\linewidth]{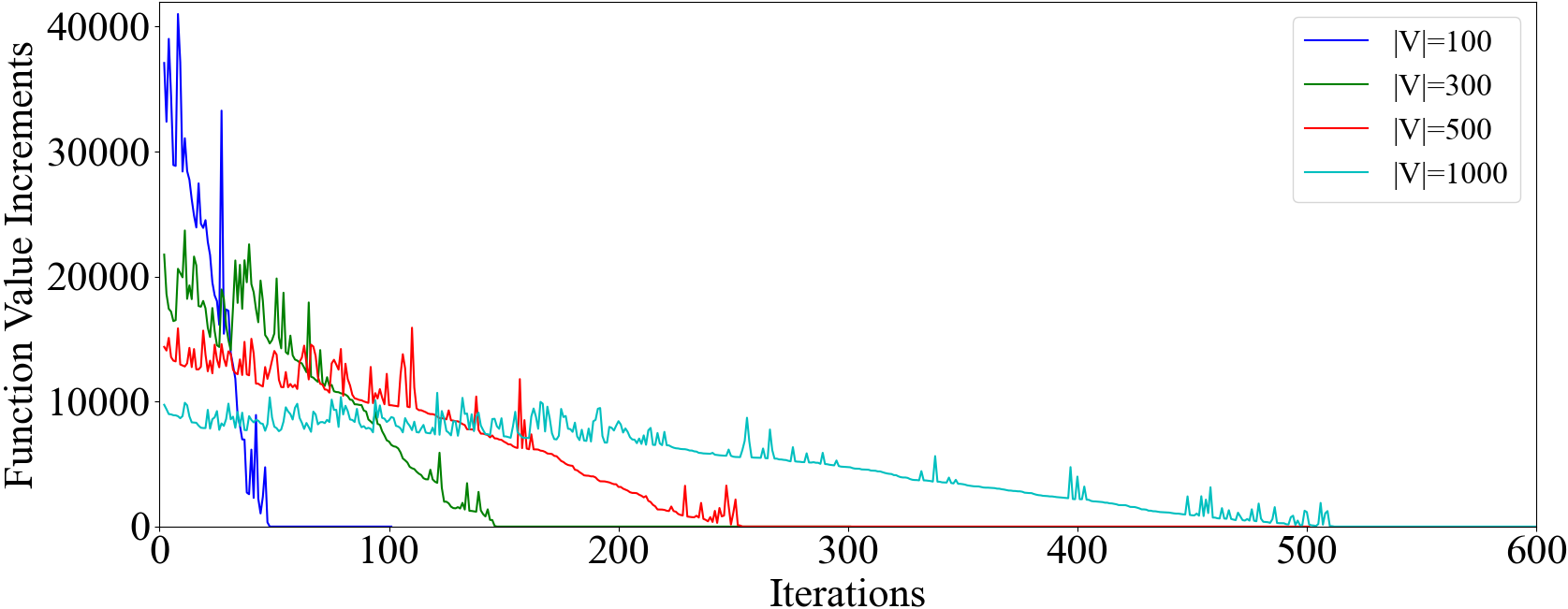}
    \caption{Convergence of the objective function with increasing iterations.}
    \label{fig:convergence}
\end{figure}

Figure~\ref{fig:convergence} presents the convergence behavior of the algorithm, tested with varying numbers of sample points: 100, 300, 500, and 1000. The results show that although increasing the number of sample points requires more iterations to reach convergence, the algorithm consistently converges within the limits of the sample size. This indicates that, even with higher sample densities, the algorithm remains robust and converges reliably, albeit with increased computational effort.

%%%%%%%%%%%%%%%%%%%%%%%%%%%%%%%%%
\section{Ablation Study on flow field}

To validate the effect of edge weighting based on the flow field on heat exchange performance, we perform an ablation study comparing the structures generated with and without edge weighting during graph initialization. This comparison evaluates the heat exchange performance between surfaces generated with and without edge weighting.

\begin{figure}[htb]
    \centering
    \includegraphics[width=\linewidth]{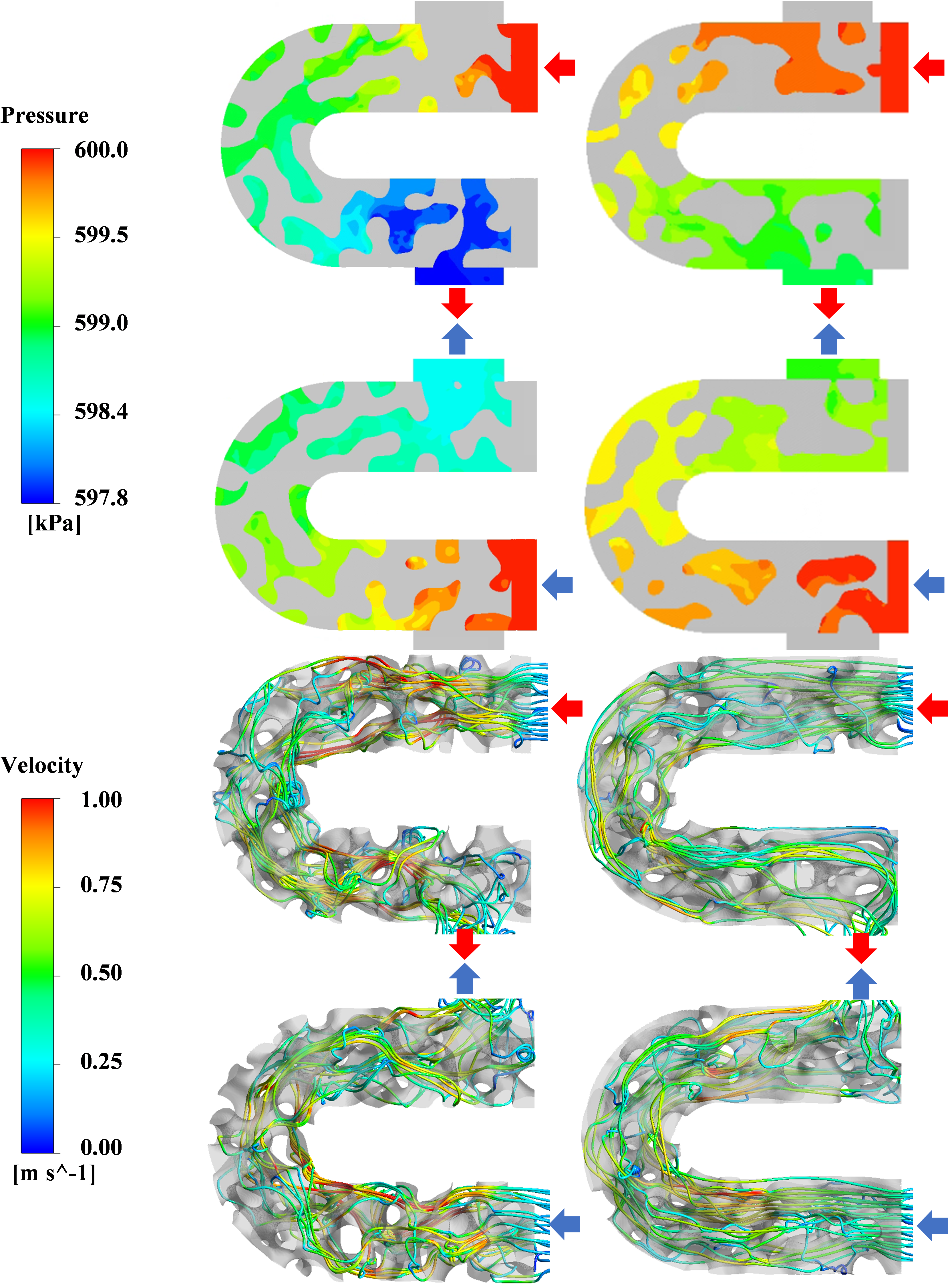}
    \caption{Ablation Study on flow field. Comparison of pressure drop and streamlines between generate structure without flow field edge weighting (left) and with flow field edge weighting $a=5$ (right), showing differences in pressure drop and streamlines. In each group, top row is hot flow channel, bottom row is cold flow channel.}
    \label{fig:ablation-flow}
\end{figure}

Figure~\ref{fig:ablation-flow} demonstrates that omitting edge weighting leads to a significant increase in pressure drop. Specifically, the average pressure drop for the structure with flow field edge weighting is 1053 Pa, compared to 1871 Pa for the structure without weighting, indicating a substantial reduction in flow resistance when the weighting is applied.

In terms of temperature difference, the structures generated with and without edge weighting have similar thermal performances. The average temperature difference for the weighted structure is 9.375 K, while the unweighted structure has a slightly higher temperature difference of 9.9585 K. This suggests that although the weighted structure has a marginally lower temperature difference, the main advantage lies in the significant reduction in pressure drop, which improves flow dynamics and energy efficiency.

%%%%%%%%%%%%%%%%%%%%%%%%%%%%%%%%%
\section{Simulation}
\label{sm:setup}

\subsection{Governing Equations}

The continuity equation, momentum equation and energy equation solved for this problem are shown in Eqs.~(\ref{eq:continuity}), (\ref{eq:momenm}), and (\ref{eq:energy}), respectively.

\begin{equation}
\label{eq:continuity}
\frac{\partial \rho}{\partial t}+\frac{\partial}{\partial x_j}\left[\rho u_j\right]=0
\end{equation}

\begin{equation}
\label{eq:momenm}
\frac{\partial}{\partial t}\left(\rho u_i\right)+\frac{\partial}{\partial x_j}\left[\rho u_i u_j+P \delta_{i j}-\tau_{j i}\right]=0, \quad i=1,2,3 \\
\end{equation}

\begin{equation}
\label{eq:energy}
\frac{\partial}{\partial t}\left(\rho e_0\right)+\frac{\partial}{\partial x_j}\left[\rho u_j e_0+u_j P-u_i \tau_{i j}\right]=0
\end{equation}

\subsection{Material Properties}
Table~\ref{tab:Materials-Properties} summarizes the material properties of the lubricating oil, fuel and Aluminum utilized in the heat exchange experiment. In this setup, the lubricating oil (lo) acts as the high-temperature fluid, while the fuel serves as the low-temperature fluid. Both fluids are subjected to identical inlet and outlet boundary conditions. The heat exchanger walls are constructed from aluminum, ensuring efficient thermal conductivity and facilitating effective heat transfer between the two fluids.

\begin{table}[htpb]
    \centering
    \caption{Dimensional physical parameters.}
    \begin{tabular}{lccc}
        \toprule
        MP\textsuperscript{a} & lo & fuel  & Aluminum\\
        \midrule
        $T_{in}$ [K] & 403.15 & 313.15 & - \\
        $P_{in}$ [MPa] & 0.6 & 0.6& - \\
        $\dot{m}_{out}$ [kg/s] & 0.13 & 0.13 & - \\
        $\rho_0$ [kg/m$^3$] & 915.74 & 801.60 & 2700.00 \\
        $\mu_0$ [Pa $\cdot$ s] & \num{2.51e-3} & \num{1.58e-3} & - \\
        $C_{p0}$ [J/(kg $\cdot$ K)] & 2102.22 & 2073.29 & 900.00 \\
        $\kappa_0$ [W/(m $\cdot$ K)] & 0.14 & 0.13 & 237.00 \\
        $Re$ [-] & 1482 & 2360 & - \\
        $L$ [m] & 0.038 & 0.038 & - \\
        \bottomrule
    \end{tabular}
    \\[0.5ex]
    \textsuperscript{a}MP: Materials Properties.
    \label{tab:Materials-Properties}
\end{table}

\subsection{Simulation Experimental Data}

Table~\ref{tab:SimulationData} presents the experimental results for the U-shape structure generated by our approach with varying numbers of sample points (300, 400, 500, 600, and 800), alongside the Gyroid-filled U-shape under identical operating conditions. The results analyze the changes in temperature ($\Delta T$) and pressure drop ($\Delta P$) of fuel and lubricating oil across different volume fractions. Our design aims to minimize $\Delta P$ while maintaining a consistent temperature difference $\Delta T$.

\begin{table}[htpb]
    \centering
    \caption{The result of comparison with uniform Gyroid structures.}
    \begin{tabular}{lccccc}
        \toprule
         & $\Delta T_{\text{lo}}$  [\wz{K}]& $\Delta T_{\text{fuel}}$ [\wz{K}]&$\Delta P_{\text{lo}}$ [\wz{Pa}]&
         $\Delta P_{\text{fuel}}$ [\wz{Pa}] &\\
        \midrule
        $\text{Gyroid}_1$ &8.11&8.138&1415&1837\\
        $\text{Ours}_1$ &7.578&7.573&912&976\\
        $\text{Gyroid}_2$ &9.427&9.321&1878& 1842\\
        $\text{Ours}_2$ &9.468&9.282&1142&964\\
        $\text{Gyroid}_3$ &10.684&10.171&2026&2270\\
        $\text{Ours}_3$ &11.287&9.071&1274&1503\\
        $\text{Gyroid}_4$ &10.766&10.676&2015&1959\\
        $\text{Ours}_4$ &11.885&11.75&1525&1530\\
        $\text{Gyroid}_5$ &12.569&12.353&2563&2091\\
        $\text{Ours}_5$ &11.644&11.427&1253&1407\\
        \bottomrule
    \end{tabular}
    \label{tab:SimulationData}
\end{table}

\section{Fabrication results}

\begin{figure}[htpb]
    \centering
    \includegraphics[width=\linewidth]{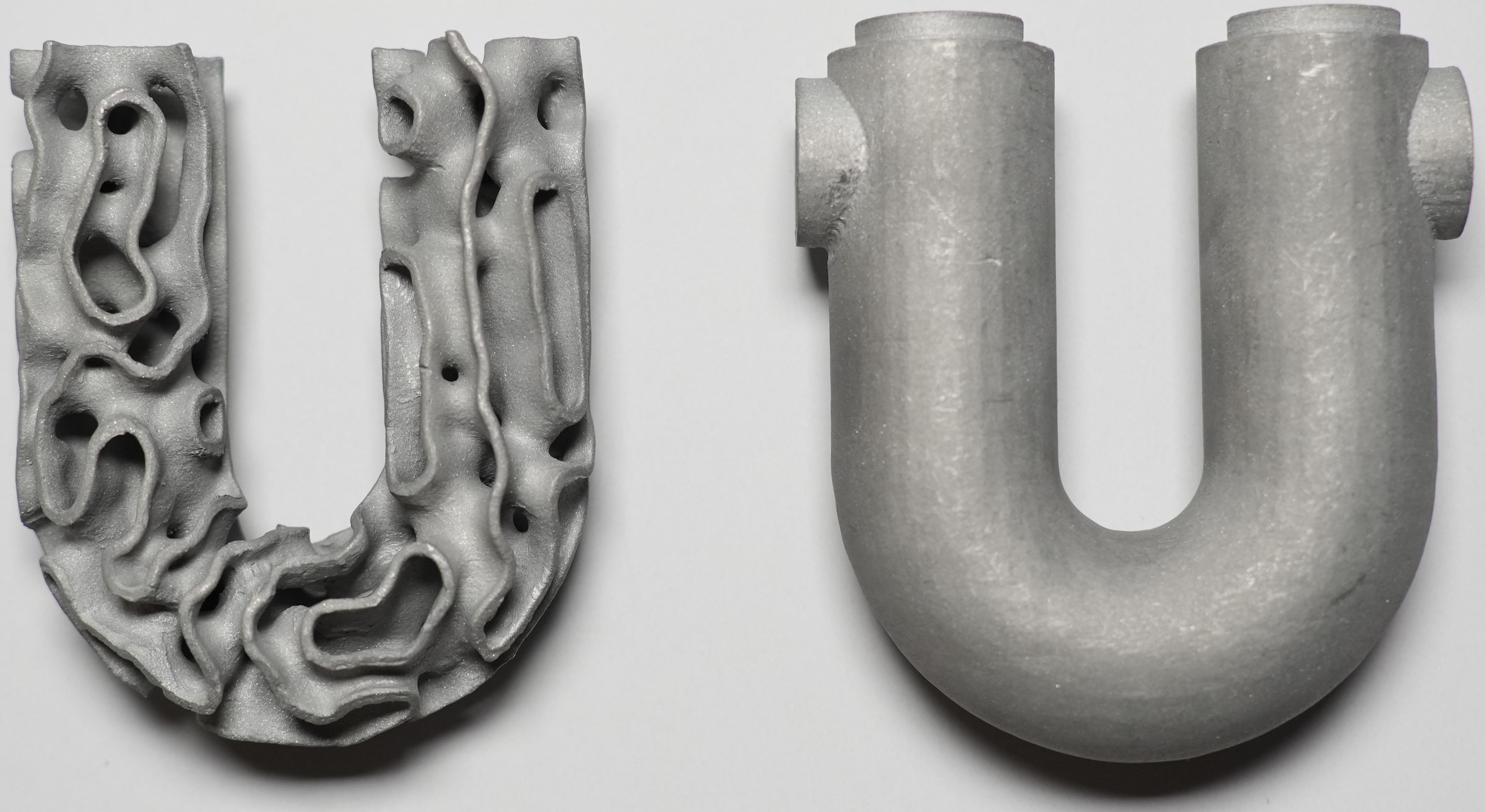}
    \caption{Metal-printed U-shape heat exchanger structure. The structure is printed without the boundary (left) and with the boundary included (right).}
    \label{fig:uShapeFab}
\end{figure}

Our method is suitable for additive manufacturing because the generated surfaces have no isolated islands or enclosed cavities, avoiding trapped powder issues. 
To validate the manufacturability of the generated structures,  we employed selective laser melting (SLM) metal additive manufacturing to fabricate a U-shape sample with 400 points, as shown in Figure~\ref{fig:uShapeFab}. 
The printed structure replicates the digital design, confirming the feasibility of our approach.

\printnomenclature